\documentclass[dvips,preprint]{imsart}
\usepackage{amsfonts}
\usepackage{graphicx}

\usepackage[english]{babel}
\usepackage[latin1]{inputenc}
\usepackage{url}
\usepackage{verbatim}
\usepackage{multirow}
\usepackage{fancyhdr}
\usepackage{lastpage}
\usepackage{latexsym}
\usepackage{amssymb}
\usepackage{longtable}
\usepackage{array}
\usepackage[dvips]{epsfig}
\usepackage{lscape}
\usepackage{amsmath}
\usepackage{epsf,epsfig,graphics,subfigure}
\usepackage{dsfont}

\def\Dy#1{\Frac{\partial #1}{\partial y}}

\def\Dy_1y_1#1{\Frac{\partial^2 #1}{\partial y_1^2}}

\def\reff#1{{\rm(\ref{#1})}}

\def\be{\begin{eqnarray}}
\def\ee{\end{eqnarray}}

\def\b*{\begin{eqnarray*}}
\def\e*{\end{eqnarray*}}

\newtheorem{Theorem}{Theorem}[part]
\newtheorem{Definition}{Definition}[part]
\newtheorem{Proposition}{Proposition}[part]

\newtheorem{Lemma}{Lemma}[part]
\newtheorem{Corollary}{Corollary}[part]
\newtheorem{Remark}{Remark}[part]
\newtheorem{Example}{Example}[part]

\makeatletter \@addtoreset{equation}{section}

\@addtoreset{Definition}{section}

\@addtoreset{Theorem}{section}

\@addtoreset{Proposition}{section}

\@addtoreset{Assumption}{section}

\@addtoreset{Corollary}{section}

\@addtoreset{Lemma}{section}

\@addtoreset{Remark}{section}

\@addtoreset{Example}{section}

\makeatother \makeatletter

\def \Int{\displaystyle\int}
\def \Frac{\displaystyle\frac}
\def \Inf{\displaystyle\inf}
\def \Sup{\displaystyle\sup}
\def \Lim{\displaystyle\lim}
\def \Liminf{\displaystyle\liminf}
\def \Limsup{\displaystyle\limsup}

\def \be{\begin{eqnarray}}
\def \ee{\end{eqnarray}}
\def \b*{\begin{eqnarray*}}
\def \e*{\end{eqnarray*}}

\def \E{\mathbb{E}}
\def \F{\mathbb{F}}

\def \N{\mathbb{N}}

\def \R{\mathbb{R}}

\def \[{[\,\!\![}
\def \]{]\,\!\!]}

\def \1{{\bf 1}}

\def \ep{\hbox{ }\hfill$\Box$}

\def\reff#1{{\rm(\ref{#1})}}

\def\Ac{{\cal A}}

\def\Cc{{\cal C}}
\def\Dc{{\cal D}}
\def\Ec{{\cal E}}
\def\Fc{{\cal F}}

\def\Hc{{\cal H}}

\def\Lc{{\cal L}}

\def\Nc{{\cal N}}

\def\Lc{{\cal L}}
\def\Pc{{\cal P}}
\def\Sc{{\cal S}}

\def\Xc{{\cal X}}

\def\Qc{{\cal Q}}

\def\Rc{{\cal R}}

 \def\N{\mathbb{N}}

\begin{document}
\begin{frontmatter}
\title{Robust utility maximization under convex portfolio constraints}
\date{}
\runtitle{}

\author{\fnms{Anis}
 \snm{MATOUSSI}\corref{}\ead[label=e1]{anis.matoussi@univ-lemans.fr}}
 \thankstext{t3}{Research partly supported by the Chair {\it Financial Risks} of the {\it Risk Foundation} sponsored by Soci\'et\'e G\'en\'erale, the Chair {\it Derivatives of the Future} sponsored by the {F\'ed\'eration Bancaire Fran\c{c}aise}, and the Chair {\it Finance and Sustainable Development} sponsored by EDF and Calyon }
\address{
 Universit\'e du Maine \\
Institut du Risque et de l'Assurance du Mans\\
Laboratoire Manceau de Math\'ematiques\\\printead{e1}
 }

 \author{\fnms{Hanen}
 \snm{Mezghani}\corref{}\ead[label=e2]{hanen.mezghani@lamsin.rnu.tn}}
%\thankstext{T1}{The work of the first author is supported by the chair \textit{risque de cr\'edit}, F\'ed\'eration bancaire Fran\c{c}aise}
\address{
University of Tunis El Manar \\
Laboratoire de Mod\'elisation
 Math\'ematique et Num\'erique\\
 dans les Sciences  de l'Ing\'enieur, ENIT\\\printead{e2}}

\author{\fnms{Mohamed}
 \snm{MNIF}\corref{}\ead[label=e3]{mohamed.mnif@enit.rnu.tn}}
\thankstext{t4}{This work was partially supported by the research project MATPYL
of the F\'ed\'eration de Math\'ematiques des Pays de la Loire }
\address{
University of Tunis El Manar \\
Laboratoire de Mod\'elisation
 Math\'ematique et Num\'erique\\
 dans les Sciences  de l'Ing\'enieur, ENIT\\\printead{e3}}
% \affiliation{Some University}

\runauthor{  A. Matoussi, H. Mezghani, M. Mnif}
\vspace{3mm}
%\title{ Maximization of recursive utilities under convex portfolio constraints}
%
%
%
%\date{}
%\begin{document}
%\author{Anis MATOUSSI,\hspace{2cm}Hanen MEZGHANI,\hspace{2cm}Mohamed MNIF}
%\maketitle

\begin{abstract}: 
We study a robust utility maximization problem from terminal
wealth and consumption  under a convex constraints on the portfolio. 
We state the existence and the uniqueness of
the consumption-investment strategy by studying the associated quadratic backward stochastic differential
equation (BSDE in short). We characterize the optimal control by using the duality method and  deriving  a dynamic
maximum principle. 
\end{abstract}

\vspace{7mm}

\noindent {\bf Key words~:} Utility maximization,  Backward Stochastic 
Differential Equations, recursive utility, model uncertainty, robust control,  
maximum principle, Forward-Backward System.

\vspace{5mm}

\noindent {\bf MSC Classification (2000)~:}
92E20, 60J60, 35B50.

\end{frontmatter}

\newpage

\section{Introduction}
The utility maximization is a basic problem in mathematical finance. It was introduced
by Merton \cite{mer}. Using stochastic control methods, he exhibits a closed
formula for the value function and the optimal proportion-portfolio when the
risky assets follow a geometric Brownian motion and the utility function is of
CRRA type.\\
In the literature, many works assume that the underlying model is exactly known.
In this paper we consider a problem of utility maximization under uncertainty.
The objective of the investor is to determine
the optimal consumption-investment strategy when the model is not exactly known.
Such problem is known as the robust utility maximization and is formulated as
\begin{eqnarray}\label{pbglobal}
\mbox{find}\,\,\Sup_{\pi}\Inf_{Q}U(\pi,Q)
\end{eqnarray}
where $U(\pi,Q)$ is the $Q$-expected utility.
The investor has to solve a sup inf problem. He considers the worst scenario by minimizing
over a set of probability measures and then he maximizes his utility.
In the literature there are two approaches to solve the robust utility maximization problems.
The first one relies on duality methods such as Quenez \cite{Quen} or Shied and Wu \cite{schi}.
They considered a set of probability measures called priors and they minimized over this set.
The second approach, which is followed in this paper, is based on the penalization method
and the minimization is taken
over all possible models such as in Anderson, Hansen and Sargent \cite{and}.
Moreover Skiadas \cite{sk03} followed the same point of view and he gave the dynamics of the control problem
via BSDE in the Markovian context. 
In our case, the $Q$-expected utility is the sum of a classical utility function and a
penalization term based on a relative entropy.
In Bordigoni et al. \cite{brom1}, they proved the existence of a unique $Q^*$ optimal model which minimizes
our cost function. They used the stochastic control techniques to study the dynamic value of the minimization problem. 
In the case of continuous filtration, they showed that the value function is the unique solution 
of a generalized BSDE with a quadratic driver.\\
In Faidi, Matoussi and Mnif \cite{fai}, they studied the maximization part of the problem \reff{pbglobal} 
in a complete market by using the BSDE approach 
as in Duffie and Skiadas \cite{duf1} and  El Karoui et al. \cite{elk}.\\
In our paper, we assume that the portfolio is constrained to take values in a given 
closed convex non-empty subset $K$ of $\R^d$. Such problem was studied when the underlying model is known by 
Karatzas, Lehoczky, Shreve and Xu \cite{kar2} in the incomplete market case and then by 
Cvitanic and Karatzas \cite{cvi} for convex constraints on the portfolio. 
Skiadas and Schroder \cite{sch03} studied the lifetime consumption-portfolio recursive utility problem under convex trading constraints. 
They used the utility gradient approach. They derived a first order conditions of optimality which take the form of a constrained 
forward backward stochastic differential equation. Wealth was computed in a
recursion starting with a time-zero value forward in time, while utility was computed in a recursion starting with a terminal date value backward
in time. 
In our context, we study the robust formulation of the consumption-investment utility problem under convex constraints on the portfolio. 
 %Moreover,  we assume that the wealth process satisfies an arbitrary uniform lower bound and the terminal wealth in non negative. 
Using change of measures and optional decomposition under constraints, we give a dual characterization of the admissible 
consumption investment strategy, then we state an existence result 
to the optimization problem where the criterion is the solution at time 0 of a quadratic BSDE with unbounded terminal condition.
To describe the structure of the solution, we use duality arguments.
The heart of the dual approach in the classical setting, when the criterion is taken under the historical 
probability measure, is to find a saddle point for the Lagrangian and apply a mini-max theorem in the infinite dimensional case.
It is appropriate to use the conjugate function of $U$ and $\bar U$. In our case, the criterion is taken under the probability measure modeling the worst scenario and the conjugate function does not appear naturally. We use the duality arguments in a different way. 
We prove the existence of a probability measure under which the budget constraint is satisfied with equality. 
%We prove the existence of a saddle point.  
Then, we derive a maximum principle which gives a necessary and sufficient conditions of optimality. 
Thanks to this result, we 
give an  implicit expression of the optimal terminal wealth and the optimal consumption rate. 
This later result is a generalization of Cvitanic and Karatzas \cite{cvi} work. \\
%Moreover, our work is a generalization  of results in
%Skiadas and Schroder \cite{sch99,sch03} and El Karoui et al. \cite{elk}.
%to the semimartingale context with unbounded terminal condition for the associated BSDE.\\
% and the method is based on the study of the associated generalized BSDE with unbounded terminal condition.\\ 
The paper is organized as follows. Section 2 describes the
model and the stochastic control problem. 
Section 3 is devoted to the existence and the uniqueness of an optimal strategy.
In section 4, we characterize the optimal consumption strategy and the optimal terminal
wealth by using duality techniques. In section 5, we relate the optimal  control  
to the solution of a forward-backward system and we study some examples. 
\section{Problem formulation}
We consider a probability space $(\Omega,\Fc,P)$ supporting a d-dimensional standard Brownian motion
$W=(W^1,...,W^d)$, over the finite time horizon $[0,T]$. We shall denote by $\F$ 
the $P$-augmentation of the filtration $(\Fc_t)_{0\leq t\leq T}=(\sigma(W_s,\,\,0\leq s\leq t))_{0\leq t\leq T}$ generated by $W$. 
We also assume that $\Fc=\Fc_T$.\\
For any probability measure $Q\ll P$
on $\Fc_T$, the density process of $Q$ with respect to $P$ is the continuous $P$-martingale
$Z^Q=(Z^Q_t)_{0\leq t\leq T}$ with
\begin{eqnarray*}
Z_t^Q=\frac{dQ}{dP}\Big |_{\Fc_t}=E_P\Big[ \frac{dQ}{dP}\Big |\Fc_t\Big].
\end{eqnarray*}
Bordigoni et al. \cite{brom1} studied a robust control problem with a
dynamic value process of the form
\begin{eqnarray}\label{robust}
Y_t= \mbox{ess}\Inf_{Q\in \Qc_f}\Big(
\frac{1}{S_t^\delta}
E_Q\Big[ \int_{t}^{T}\alpha S_s^\delta \check U_sds+ \bar \alpha S_T^\delta \bar U_T \Big |\Fc_t \Big]
+\beta E_Q\Big[\Rc_{t,T}^{\delta}(Q)
|\Fc_t \Big]
\Big),
\end{eqnarray}
where
\begin{eqnarray*}
\Qc_f=\{Q|Q\ll P, Q=P \mbox{ on }\Fc_0 \mbox{ and } H(Q|P):=E_Q[\log \frac{dQ}{dP}]<\infty\},
\end{eqnarray*}
$\alpha$ and $\bar \alpha$ are non negative constants, $\beta\in (0,\infty)$,
$\delta= (\delta_t)_{0\leq t\leq T}$ and
$\check U= (\check U_t)_{0\leq t\leq T}$ are $\F$ progressively measurable processes, $\bar U_T $
is a $\Fc_T$ measurable random variable,
$S_t^\delta=e^{-\int_0^t\delta_sds}$ is the discounting factor and $\Rc_{t,T}^{\delta}$ is the penalization term
which is the sum of the entropy rate and the terminal entropy:
\begin{eqnarray*}
\Rc_{t,T}^{\delta}=\frac{1}{S_t^\delta}\int_{t}^{T}\delta_s S_s^\delta \log \frac{Z_s^Q}{Z_t^Q} ds+
\frac{S_T^\delta}{S_t^\delta}\log \frac{Z_T^Q}{Z_t^Q}.
\end{eqnarray*}
We define the following spaces:\\
$L_+^0(\Fc_T)$ is the set of nonnegative $\Fc_T$ measurable random variables.\\[0.2cm]
$L^{\exp} $ is the space of all ${\cal{ F}}_T$measurable random variables $X $ with
 $$E_{P}\left[\exp\left(\gamma \vert X  \vert\right)\right]<\infty \qquad \hbox{ for all }
 \gamma>0.$$
 $D^{\exp}_0$ is the space of all progressively measurable processes $X={(X_t)}_{0\le t\le T}$ with
 $$E_{P}\left[\exp\left(\gamma {~ \rm{ess}\sup}_{0\le t\le T}|X_t|\right)\right]<\infty \qquad
 \hbox{ for all } \gamma>0.$$
 $D^{\exp}_1$ is the space of all progressively measurable processes $X={(X_t)}_{0\le t\le T}$ such that
$$E_{P}\left[\int_0^T\exp\left(\gamma |X_s|\right)ds\right]<\infty \qquad \hbox{ for all } \gamma >0.$$
$H^2_T(\R^d)$ is the set of progressively measurable processes $\R^d$ valued $Z={(Z_t)}_{0\le t\le T}$ such that
$$ ||Z||_{H^2}^2:=E_{P}\big[ \int_0^T |Z_t|^2dt\big]<\infty. $$
We shall assume the boundedness on the discounting factor and the exponential integrability of the utility processes i.e.\\
{\bf (H1)} $0\leq \delta\leq ||\delta ||_\infty$ for some constant $||\delta ||_\infty$.\\
{\bf (H2)} $\check U \in  D_1^{\exp}$ and $\bar U_T \in  L^{\exp}$.\\
Under the boundedness on the discounting factor {\bf (H1)} and  the exponential integrability of the utility processes {\bf (H2)}, Bordigoni et al. \cite{brom1} (Theorem 6,Theorem 12 and Proposition 16) proved the existence and the uniqueness  of an optimal
probability measure $Q^*$ of  the problem \reff{robust}. They showed that the dynamics of $(Y_t)_{t\in [0,T]}$ satisfies the following BSDE
\begin{eqnarray}
dY_t &=&(\delta_t Y_t -\alpha \check U_t)dt +\frac{1}{2\beta}|Z_t^Y|^2dt+Z_t^{Y \prime}dW_t,\label{BSDE}\\
Y_T&=&\bar \alpha \bar U_T,\label{BSDEct}
\end{eqnarray}
where $|.|$ stands the euclidean norm and the notation $^\prime$ denotes the transposition operator. They  established for $Y$ the recursive relation
\begin{eqnarray}\label{recurciverelation}
Y_t=-\beta \log E_P\Big[\exp\Big(\frac{1}{\beta} \int_t^T (\delta_sY_s-\alpha \check U_s)ds
- \frac{1}{\beta}\bar \alpha \bar U_T \Big)\Big|\Fc_t\Big].
\end{eqnarray}
They proved that there exists a unique pair $(Y,Z^Y)\in D_0^{exp}\times H^2_T(\R^d)$ that solves \reff{BSDE}-\reff{BSDEct}.\\
Moreover, they showed that the density of the probability measure $Q^*$ is a true martingale and is given by
\begin{eqnarray}\label{minimizingmeasure}
Z_t^{*}=\Ec_t(-\frac{1}{\beta}M^Y),0\leq t\leq T,
\end{eqnarray}
where $M^Y_t=\int_0^tZ_s^{Y \prime}dW_s$; $t\in [0,T]\, dt\otimes dP $ and  $\Ec$ denotes the stochastic exponential.\\
From now, we are interested in the problem of utility maximization.
Let us consider an investor who can consume between time $0$ and time $T$.
We denote by $c=(c_t)_{0\leq t\leq T}$ the consumption rate.
We consider a financial market consisting of a bond and $d$ risky assets. Without loss of generality, we assume that the bond is constant.
The risky assets $S:=(S^1,...,S^d)$ evolve according to the stochastic differential equations
\begin{eqnarray*}
dS^i_t&=&S^i_t\big(b^i_tdt+\sum_{j=1}^{d}\sigma^{ij}_tdW^j_t\big),\,\,\,S^i_0=1,\,i=1...d.
\end{eqnarray*}
We assume that the process $b=(b^1_t,...,b^d_t)_{t\in [0,T]}$ (vector of instantaneous yield) and the process $\sigma=\Big((\sigma_t^{ij})_{1\leq i,j\leq d}\Big)_{t\in [0,T]}$ (volatility matrix) are $\F$ progressively measurable. 
%and satisfy
%\begin{eqnarray}
%\label{elliptic}
%\sigma_t^{-1}\xi &\geq &\frac{1}{\epsilon} \xi\,\,\,dt\otimes dP \,a.s. \forall (t,\xi)\in [0,T]\times \R^d,
%\end{eqnarray}
%for a given real constant $\epsilon >0$, as well as
%\begin{eqnarray*}
%E\int_0^Tr_tdt < \infty.
%\end{eqnarray*}
We shall assume throughout that the relative risk process
\begin{eqnarray*}
\theta_t:=\sigma^{-1}_tb_t,
\end{eqnarray*}
satisfies the integrability condition
\begin{eqnarray*}
\int_0^T||\theta_t||^2dt < \infty,\, P\,a.s.
\end{eqnarray*}
We denote by
$H=((H^1_t,...,H^n_t)_{t\in[0,T]})^{\prime}$ the investment strategy representing the amount of each
asset invested in the portfolio.
We shall fix throughout a nonempty, closed, convex set $K$ in $\R^d$ containing 0, and denote by
\begin{eqnarray*}
\delta^{supp}(x):=\delta^{supp} (x|K):=\Sup_{y\in K}(-y^{\prime}x):\R^d\longrightarrow \R \cup \{+\infty\}
\end{eqnarray*}
the support function of the convex set $-K$. This is a closed, positively homogeneous, 
proper convex function on $\R^d$  finite on its effective domain (Rockafellar  \cite{R:book} p. 114) 
\begin{eqnarray*}
\tilde K&:=&\{ x\in \R^d,\,\,\delta^{supp} (x|K)<\infty\}\\
&=&\{ x\in \R^d,\,\,\mbox{ there exists }\beta\in\R \mbox{ s.t. }-y^{\prime}x\leq \beta,\,\, \forall y\in K\},
\end{eqnarray*} 
which is a convex cone (called the barrier cone of $-K$).\\
We assume that 
\begin{eqnarray}\label{fonctsupp}
\tilde{K}\,\, \mbox{is closed and the function } \delta^{supp}(.|K) \mbox{ is continuous on }\tilde K.
\end{eqnarray}
{\bf Examples}
\begin{itemize}
\item \underline{$K$ is a linear space: unconstrained portfolio}\\
$K=\R^d$
\begin{eqnarray*}
\left\{
\begin{array}{ll}
\delta^{supp}(x)&=0\mbox{ if },x=0\\
\delta^{supp}(x)&=\infty\mbox{ otherwise },
\end{array}
\right.
\end{eqnarray*} 
and $\tilde K=\{ 0\}$.\\
\item \underline{$K$ is a convex closed cone in $\R^d$: short selling contract}\\
$K=\{\pi \in \R^d;\,\pi_i\geq 0\,\,\forall \,i=1,\dots,d\}$.\\
\begin{eqnarray*}
\left\{
\begin{array}{ll}
\delta^{supp}(x)&=0\mbox{ if },x\in K\\
\delta^{supp}(x)&=\infty\mbox{ otherwise },
\end{array}
\right.
\end{eqnarray*} 
and $\tilde K=K.$\\
\item \underline{$K$ is a convex closed set in $\R^d$: rectangular constraints}\\
$K= \displaystyle{\prod_{i=1}^{d}} K_i$ where $K_i=[\alpha_i,\beta_i]$; $ -\infty<\alpha_i\leq 0 \leq \beta_i\leq +\infty$\\
$\delta^{supp}(x)=\sum_{i=1}^{d}\beta_ix_i^{-}-\sum_{i=1}^{d}\alpha_ix_i^{+}$ and $\tilde K=\R^d.$\\
\end{itemize}
In all these examples $\tilde K$ is closed and the support function is continuous on $ \tilde K$ i.e. Assumption (\ref{fonctsupp}) is satisfied. \\\\ 
The investment strategy is constrained to remain in the convex set K.
We denote by $\tilde \Cc$ and $\tilde \Hc$ the following sets
\begin{eqnarray*}
\tilde \Cc&:=&\{ c=(c_t)_{t\in [0,T]}\,\F-\mbox{progressively measurable },\,\, c_t \geq 0\,dt\otimes dP \mbox{ a.e. and } \int_0^T c_tdt <\infty\},\\
\tilde \Hc&:=&\{ H=(H_t)_{t\in [0,T]}\,\F-\mbox{progressively measurable}, \,\R^d\mbox{ valued and} \,\, H^{\prime}\mbox{diag}(S)^{-1}\in L(S)\\
& &\mbox{ and } H_t\in K\, dt\otimes dP \mbox{ a.e.}\},
\end{eqnarray*}
where $L(S)$ denotes the set of $\F-$ progressively measurable processes, $\R^d$ valued such that the stochastic integral with respect to $S$ is well-defined.\\ 
Given an initial wealth $x\geq 0$ and a policy $(c,H)\in \tilde \Cc\times \tilde \Hc$, the wealth process at time $t$ follows the dynamics given by:
\begin{eqnarray}
d X_t^{x,c,H}= H_t^{\prime}\mbox{diag}(S_t)^{-1}dS_t -c_tdt,\,\,\,\,\, X_0^{x,c,H}=x.
\end{eqnarray}
%We impose the following constraint on the wealth process
%\begin{eqnarray}
% X_t^{x,c,H}&\geq& d \mbox{ a.s. } \forall t\in [0,T) \mbox{ for some }d\in \R,\label{constraint1}\\
% X_T^{x,c,H}&\geq& 0 \mbox{ a.s. }\label{constraint2}
%\end{eqnarray}
The investor has preferences modeled by the utility functions $U$ and $\bar U$. \\
We define the set of admissible strategies as follows: 
\begin{Definition}
(i) We denote by $\Ac$ as the set of all processes  $(c,\xi)\in \tilde \Cc\times L^0_+( \Fc_T)$
such that  
%constraints \reff{constraint1}-\reff{constraint2} hold as well as 
the families
\begin{eqnarray}
& &\{ \int_0^T\exp{(\gamma  |U(c_s)|)ds} \, : \, c\in \tilde \Cc \}=:\Cc^{ad}\label{unifinteg1}, 
%& &\{ \int_0^T\exp{(\gamma  |U^{'}(c_s)|)ds} \, : \, c\in \tilde \Cc  \}\label{unifinteg2},
\end{eqnarray}
\begin{eqnarray}
& &\{ \exp{(\gamma |\bar U(\xi)|)} \, : \, \xi\in L^0_+(\Fc_T)\}=:\Hc^{ad} ,\label{unifinteg3} 
%& &\{ \exp{(\gamma |\bar U^{'}(\xi)|)} \, : \, \xi\in L^0_+(\Fc_T) \}\label{unifinteg4}
\end{eqnarray}
are uniformly integrable for all $\gamma >0$.\\
(ii)Given an initial wealth $x\geq 0$, we define the set $\Ac(x)$ as the set of all processes  $(c,\xi)\in \Ac$
such that there exists
$H\in \tilde \Hc$ satisfying $X_T^{x,c,H}\geq \xi$. 
%(ii) $\Hc $ consists of all processes $H\in \tilde \Hc$ such that \reff{unifinteg3} and \reff{unifinteg4} are checked.
\end{Definition}
We shall assume\\
{\bf (H3)} The set of controls $\Ac$ is non-empty.
%\begin{Remark}
%The set $\Cc\times \Lc$ is a convex one when the utilities are power functions.
%\end{Remark}
\begin{Remark}
We will see in the examples that Assumption {\bf (H3)} is satisfied when the utility function $\bar U$ is logarithmic.
\end{Remark}
%\begin{Remark}
%Hu el al. (2005) assume similar condition of uniform integrability. In the literature, this condition coincides with the notion of class $D$ (see Dellacherie and Meyer (1975)).
%\end{Remark}
%\begin{Remark}
%We can weaken the assumption of the uniform integrability of the families \reff{unifinteg2} and \reff{unifinteg4}. In fact, We only need that $\int_0^T |U^%{'}(c_t)|dt \in L^4(P)$ for all $c\in \tilde \Cc$ and $\bar U^{'}(X_T^{x,c,H})\in L^4(P)$ for all $(c,H)\in  \tilde \Cc\times \tilde \Hc$ to prove the uppe%r-semicontinuity of the value function (See Lemma \ref{regularity}).
%\end{Remark}
The problem of optimal consumption-investment is formulated as
\begin{eqnarray}\label{reward}
V(x)=\Sup_{(c,\xi)\in  \Ac(x) }Y_0^{x,c,\xi},\,\,\,x\in \R_+,
\end{eqnarray}
where the dynamics of $Y^{x,c,\xi}=(Y_t^{x,c,\xi})_{0\leq t\leq T}$ is given by
\begin{eqnarray}
dY_t^{x,c,\xi}&=&(\delta_tY_t^{x,c,\xi} -\alpha U(c_t))dt +\frac{1}{2\beta}|Z_t^{x,c,\xi}|^2dt+Z_t^{x,c,\xi\prime}dWt\label{bsdech}\\
Y_T^{x,c,\xi}&=&\bar \alpha \bar U(\xi).\label{bsdechct}
\end{eqnarray}
The following result is a strict comparison theorem for the BSDE \reff{BSDE}-\reff{BSDEct}. We give the proof in the Appendix.
\begin{Theorem}\label{ctheorem}
We consider $(Y^1, Z^1)$ and $(Y^2, Z^2)$ two solutions of the BSDE \reff{BSDE}-\reff{BSDEct}
associated to $(\check U^1,\bar U_T^1)$ and $(\check U^2,\bar U_T^2)$ respectively.
We assume that the boundedness on the discounting factor {\bf (H1)} and the exponential integrability of the utility processes {\bf (H2)} hold and that
\begin{eqnarray}
\check U^1_t&\leq& \check U^2_t \,\,dt\otimes dP \mbox{ a.e.},\,t\in[0,T],\label{c1}\\
\bar  U_T^1&\leq& \bar U^2_T \,\,dP \mbox{ a.s. }\label{c2}
\end{eqnarray}
Then, we have
\begin{eqnarray*}
Y_t^1\leq Y_t^2 \,\,dt\otimes dP \mbox{a.e.},\,t\in[0,T].
\end{eqnarray*}
Moreover, the comparison is strict. If $Y_0^1=Y_0^2$, then $Y_t^1=Y^2_t$, $t\in [0,T]$ $dt\otimes dP$ a.e. In particular, if $P(\bar U_T^1<\bar U_T^2)>0 $ 
or $E[\int_0^T {\bf 1}_{\{\check U_t^1< \check U_t^2\}}dt]>0$, then $Y_0^1<Y_0^2$.  
\end{Theorem}
Also, we have a continuity result for the solution of the BSDE \reff{bsdech}-\reff{bsdechct} which will be useful to show the regularity of $ (c,\xi)\rightarrow Y_0^{x,c,\xi}$ and to prove the dynamic maximum principle. 
The proof is given in Faidi et al. \cite{fai} (Proposition 3.2 pp. 1024) and for  sake of completeness, we give the proof in the Appendix.
\begin{Proposition}\label{monotonieYn}
We assume the boundedness on the discounting factor {\bf (H1)}, the exponential integrability of the utility processes {\bf (H2)}
and the set of controls $\Ac$ is non-empty {\bf (H3)}.
%We assume {\bf (H1)} and {\bf (H3)}.
Let  $(c,\xi)\in  \Ac $ and $(c^n,\xi^n)_{n\in \N}\ \Ac^{\N} $.\\
(i) If $\xi^n\searrow \xi$ $dP$ a.s. and $c^n_{t}\searrow c_t $,
$0\leq t\leq T$, $dt\otimes dP$ a.e. when $n$ goes to infinity,
then $Y_t^{x, c^n,\xi^n}\searrow Y_t^{x, c,\xi}$ ,
$0\leq t\leq T$, $dt\otimes dP$ a.s. when $n$ goes to infinity.   \\
(ii) If $\xi^n\nearrow \xi$ $dP$ a.s. and $c^n_{t}\nearrow c_t$,
$0\leq t\leq T$, $dt\otimes dP$ a.e. when $n$ goes to infinity,
then $Y_t^{x, c^n,\xi^n}\nearrow Y_t^{x,c,\xi}$ ,
$0\leq t\leq T$, $dt\otimes dP$ a.s. when $n$ goes to infinity.
\end{Proposition}
\section{Optimum Strategy Plan}
In this section, we will study the existence and the unicity of an optimal consumption-investment strategy. First, we give a dual characterization of the set of admissible strategies in terms of a set of probability measures.\\   
By the martingale representation theorem for Brownian motion (see
e.g. Karatzas and Shreve \cite{karshr}), any probability
measure equivalent to $P$ has a density process in the form:
\begin{eqnarray} \label{znu}
Z^\nu &=& \Ec \left( - \int (\theta + \sigma^{-1} \nu)^{\prime}dW \right),
\end{eqnarray}    
where $\nu$ lies in the set $\Nc$ of $\R^d$-valued $\F-$ progressively measurable process such that\\ 
$\int_0^T|\sigma_t^{-1}\nu_t|^2dt$ $<$ $\infty$ P a.s. and $E[Z_T^\nu]=1$.  We define  the process  $(V_t^{c,H})_{t\geq 0}$  by $$V_t^{c,H}:=
\int_0^t H_s^{\prime} \mbox{diag}(S_s)^{-1}dS_s- \int_0^t c_sds .$$ 
By Girsanov's theorem, the Doob-Meyer decomposition of $(V_t^{c,H})_{t\geq 0}$ under $P^\nu$ $=$
$Z_T^\nu.P$,\\
for $\nu$ $\in$ $\Nc$ and 
$(c,H)\in \tilde C  \times \tilde \Hc$, is given by:
\begin{eqnarray} \label{deccons}
V^{c,H}_t &=& \int_0^t H_s^{\prime}\sigma_s dW_s^\nu - \int_0^t c_sds + A_t^{\nu,H},\,\, t\in[0,T],
\end{eqnarray}
where $W^\nu$ is a $P^\nu$-Brownian motion and the process $A^{\nu,H}=(A_t^{\nu,H})_{0\leq t\leq T}$ is given by\\
 $A_t^{\nu,H}:= \int_0^t \left( -H_s^{\prime}\nu_s \right)ds$, $0\leq t\leq T$.\\
%given by
%\begin{eqnarray*}
%A_t^{\nu,H} &=& \int_0^t \left( -H_s'\nu_s \right)ds.
%\end{eqnarray*} 
We introduce the following set of probability measures:
\begin{Definition}\label{upperbound}
(i)We denote by $\bar{\Pc}^{0}$ the class of all probability measures $P^{\nu} \sim P$ with the following property:
there exists a nondecreasing predictable process $A$  such that
\begin{eqnarray}\label{superMP} 
V^{c,H}_t + \int_0^t c_sds - A_t,\,\, t\in[0,T],
\end{eqnarray}
is a $P^\nu$-local supermartingale for any $(c,H)\in \tilde{\Cc}\times \tilde{\Hc}$.\\
(ii) The upper bound process denoted by $A_t(\nu)$ is a nondecreasing predictable process with $ A_0(\nu)=0$ which satisfies (\ref{superMP}) and 
\begin{eqnarray}\label{propupper} 
(A_t-A_t(\nu))_{t\in[0,T]}
\end{eqnarray}
is nondecreasing for all $A$ nondecreasing process satisfying (\ref{superMP}).
\end{Definition}
Therefore, by F\"ollmer and Kramkov  (\cite{ folkra}, Lemma 2.1) the probability measure  $P^{\nu}$ belongs to $\bar  \Pc^0 $ 
if and only if  there is an upper bound for all predictable processes arising in the Doob-Meyer decomposition
of the semimartingale $V^{c,H}$ under $P^\nu$, denoted in our case by $A^{\nu,H}$. In this case the upper variation is equal to this upper bound. 
Thanks  again to Lemma 2.1 in \cite{ folkra}, the set
$\bar \Pc^0$ consists of all probability measures $P^\nu$ for  $\nu \in \Nc(\tilde K)$ where  $$\Nc(\tilde K):=\{\nu \in \Nc~: \, \nu \in  \tilde K \; \mbox{ and} \;  \int_0^T \delta^{supp}(\nu_t)dt<\infty \;  \} .$$ 
The upper variation process is given by:
\begin{eqnarray*}
A_t(\nu) &=& \int_0^t\delta^{supp}(\nu_s) ds,\,\, t\in[0,T].
\end{eqnarray*}
%Moreover, for any non-decreasing predictable process $A$, such that
%\begin{eqnarray*} 
%V^{c,H}_t + \int_0^t c_sds - A_t,\,\, t\in[0,T],
%\end{eqnarray*}
%is a $P^\nu$-local super-martingale, we have
%\begin{eqnarray}\label{propupper} 
%(A_t-A_t(P^\nu))_{t\in[0,T]} \mbox{ is non-decreasing}.
%\end{eqnarray}
We fix $\eta>1$ and $\bar \eta>1$. We denote by 
\begin{eqnarray*}
G_{\mbox{equi}}:=\{g\,:\,[0,T]\rightarrow \R^d \,\mbox{s.t.}\, \delta^{supp}(g)\, \mbox{is equi-integrable with respect to the Lebesgue measure on [0,T]}\},
\end{eqnarray*}
\begin{eqnarray*}
\Theta_{ad}=\{\nu \in \Nc(\tilde K)\,\mbox{s.t.}\, \Sup_{\nu}E[(Z_T^{\nu})^{\eta}]<\infty,\,\,\Sup_{\nu}E[(Z_T^{\nu})^{1-\bar \eta}]<\infty\, \mbox{and}\, \nu\in  G_{\mbox{equi}}\,P.a.s\}.
\end{eqnarray*}
We denote by $\Pc^0$ the subset of elements $P^{\nu}\in \bar \Pc^0$:
\begin{eqnarray*} 
\Pc^0=\{P^{\nu}\in \bar \Pc^0\,\mbox{such that},\,\nu\in \Theta_{ad} \}
\end{eqnarray*}
%We fix $\eta>1$ and $\bar \eta>1$. We denote by $\Pc^0$ the subset of elements $P^{\nu}\in \bar \Pc^0$ %such that  
%$\Sup_{\nu}E[(Z_T^{\nu})^{\eta}]<\infty$  and $\Sup_{\nu}E[(Z_T^{\nu})^{1-\bar \eta}]<\infty$ 
%and the family $\{\delta^{supp}(\nu)\}_{\nu}$ is equi-integrable with respect to the Lebesgue measure on %[0,T]. 
\begin{Remark}
Such restriction is needed to characterize the optimal strategy of consumption investment (see Theorem \ref{bd}).
\end{Remark} 
As is Pham \cite{pha00b} and in order to obtain a dual characterization of dominated random variables $\xi\,\,\Fc_T$-measurable by a controlled process i.e. there exists $u_0$ and an admissible strategy of consumption-investment denoted by $(c,H)$ such that $\xi\leq X_T^{u_0,c,H}$, 
we shall assume 
\begin{eqnarray}\label{uupper1} 
E_{P^\nu}[A_T(\nu)]<\infty
\mbox{ for all } P^{\nu}\in \Pc^0,
\end{eqnarray} 
\begin{eqnarray}\label{uupper2} 
\mbox{ess}\Inf_{P^{\nu}\in \Pc^0}E_{P^\nu}[A_T(\nu)\mid \Fc_t] \mbox{ is bounded in (t,w)}. 
\end{eqnarray}
All these conditions are satisfied in the example of the last section. 
\begin{Proposition} \label{dualdomin}
We assume that the set of controls $\Ac$ is non-empty {\bf(H3)}. Let $x$ $\in$ $\R_+$ and $(c,\xi)\in \tilde{\Cc} \times L^0_+(\Fc_T)$. Then there exists $H\in \tilde \Hc$ such that $\xi\leq X_T^{x,c,H}$ if and only if 
\begin{eqnarray} \label{vhx}
v(c,\xi) := \sup_{P^{\nu}\in\Pc^0} 
E_{P^{\nu}} \left[\xi+\int_0^Tc_tdt
-A_T({\nu}) \right] \; \leq \; x.
\end{eqnarray}
\end{Proposition}
{\bf Proof.}
{\it Necessary condition.} We consider $(c,\xi)\in \tilde{\Cc} \times L^0_+(\Fc_T)$ and  $P^{\nu}$ $\in$ $\Pc^0$. There exists $H\in \tilde {\Hc}$ such that $\xi\leq X_T^{x,c,H}$. 
Since $ X_{.}^{x,c,H}+\int_0^{.}c_tdt -A_{.}(\nu)$ is a $P^{\nu}$-local supermartingale,
there exists a sequence of stopping time $(\tau_n)_{n \in \N}\nearrow \infty $ when $n$ goes to infinity, such that
\begin{eqnarray*} 
E_{P^{\nu}} \left[ X_{T\wedge \tau_n }^{x,c,H}+\int_0^{T\wedge \tau_n}c_tdt -A_{T\wedge \tau_n}(\nu) \right] \leq  x.
\end{eqnarray*}
By condition \reff{uupper1}, the nondeceasing property of $(A_t(\nu))_t$ and since $X_T^{x,c,H}+\int_0^Tc_sds$ 
is nonnegative, Fatou's lemma yields that 
\begin{eqnarray*} 
\liminf_{n \rightarrow \infty} E_{P^{\nu}} \left[ X_{T\wedge \tau_n }^{x,c,H}+\int_0^{T\wedge \tau_n}c_tdt -A_{T\wedge \tau_n}(\nu) \right] \geq  E_{P^{\nu}} \left[\liminf_{n \rightarrow \infty}\big( X_{T\wedge \tau_n }^{x,c,H}+\int_0^{T\wedge \tau_n}c_tdt -A_{T\wedge \tau_n}(\nu)\big) \right].
\end{eqnarray*}
We have $T\wedge \tau_n\nearrow T \,dP\, a.s.$ and  $A_{T\wedge \tau_n}(\nu)\nearrow  A_T(\nu)\,dP\, a.s.$, when $n$ goes to infinity. We deduce that~:
\begin{eqnarray*}  
E_{P^{\nu}} \left[ X_T^{x,c,H}+\int_0^Tc_tdt -A_T(\nu) \right] \leq  x,
\end{eqnarray*}
for all $P^{\nu}$ $\in$ $\Pc^0$. This shows that $v(c,\xi)\leq x$.\\
\vspace{1mm}
{\it Sufficient condition.} 
Consider the random variable  $g$ $=$ $\xi+\int_0^Tc_tdt$. Since 
\begin{eqnarray*}
v(c,\xi) = \sup_{P^{\nu}\in\Pc^0} E_{P^{\nu}} \left[g-A_T(\nu) \right] \; \leq \; x\;<\;\infty,
\end{eqnarray*}
then by the stochastic control Lemma A.1 of F\"ollmer and Kramkov
\cite{ folkra}, there exists a RCLL version of the process:
\begin{eqnarray}\label{eqV}
V_t &=& {\rm ess}\sup_{P^{\nu}\in\Pc^0} E_{P^{\nu}}
\left[g- A_T({\nu})+A_t({\nu}) \vert \Fc_t\right] \;\;0\leq t\leq T.
\end{eqnarray}
Moreover, for any ${P^{\nu}}\in \Pc^0$, the process $(V_t-A_t({\nu}))_{t\in [0,T]}$ is a
${P^{\nu}}$-local supermartingale. By condition \reff{uupper2}, the process V is bounded from below. 
Using the optional decomposition under constraints of F\"ollmer and Kramkov 
(see their Theorem 3.1), 
the process $V$ admits a decomposition:
\begin{eqnarray*}
V_t &=& v(c,\xi)+U_t-C_t,\,\,t\in [0,T]
\end{eqnarray*}
where $U\in \tilde \Sc:=\{X^{x,c,H}+\int cdt -x,\,c\in \tilde \Cc,\,H\in \tilde \Hc\}$ and $C$ is an (optional) nondecreasing process with 
$C_0=0$. Hence there exists $H\in  \tilde \Hc$ such that $(c,X_T^{x,c,H})\in \Ac(x)$ and
\begin{eqnarray}  \label{decomposition}
V_t &\leq& X_t^{x,c,H}+\int_0^t c_s ds,\;\;\; a.s. \;\;  0\leq t\leq T.
\end{eqnarray}
Using equation \reff{eqV} for $t=T$ and inequality \reff{decomposition}, we obtain that
\begin{eqnarray*}
V_T &:= & \xi+\int_0^Tc_tdt\leq X_T^{x,c,H}+\int_0^T c_sds ,\;\;\; 0\leq t\leq T,
\end{eqnarray*}
and so $\xi\leq X_T^{x,c,H}$ where $H\in \tilde \Hc$ and the proof is ended.
\ep\\
As an immediate consequence of the last proposition, we have a necessary and sufficient 
condition in terms of the set $\Pc^0$ for the set of admissible strategies $\Ac(x)$.
\begin{Corollary}
We assume that the set of controls $\Ac$ is non-empty {\bf(H3)}. For all $x\in \R_+$, $\Ac(x)$ is non empty if and only if 
\begin{eqnarray*}
v(0,0) = \sup_{P^{\nu}\in\Pc^0} 
E_{P^{\nu}} \left[-A_T({\nu}) \right] \; \leq \; x.
\end{eqnarray*}
\end{Corollary}
\noindent{\bf Proof.}
Suppose that $\Ac(x)$ is non empty. Then, there exists $H\in \tilde \Hc$ such that $\xi=0\leq X_T^{x,0,H}$.
By Proposition \ref{dualdomin}, we have $v(0,0)\leq x$. Conversely, suppose that $v(0,0)\leq x$. 
By Proposition \ref{dualdomin},  there exists $H\in \tilde \Hc$ such that $\xi=0\leq X_T^{x,0,H}$ and in particular $\Ac(x)$ is non empty.
\ep\\
%We say that the strategy $(c,\xi)$ is admissible if $(c,\xi)\in \Ac(x)$ and  $v(c,\xi)\leq x$. 
%We denote by $\hat \Ac(x)$ the set of admissible strategies.\\
We need the following technical lemmas related to the closeness and the convexity of the set $ \Ac(x)$.
We shall assume that the utility functions satisfy the following conditions~:\\
{\bf (H4)}(i) $U:\,\R_+\longrightarrow \R$ and $\bar U:\,\R_+\longrightarrow \R$ are $C^1$ on the sets $\{U <\infty\}$ and
$\{\bar U <\infty\}$ respectively, strictly increasing and concave.\\
(ii) $U$ and $\bar U$ satisfy  the usual Inada conditions i.e.
$U^{'}(\infty)=\bar U^{'}(\infty)=0$ and $U^{'}(0)=\bar U^{'}(0)=\infty$.\\
Also we assume  the quasi concavity of the absolute value of the utility functions i.e.\\
{\bf (H5)} For all $\lambda\in [0,1]$, $z_1\geq 0$ and $z_2\geq 0$, we have $|U(\lambda z_1+ (1-\lambda)z_2)|\leq \max(|U(z_1)|,|U(z_2)|)$ 
and $|\bar U(\lambda z_1+ (1-\lambda)z_2)|\leq \max(|\bar U(z_1)|,|\bar U(z_2)|)$.\\
%{\bf (H5)(ii)} The functions U and $\bar U$ are convex.\\
\begin{Remark}
%Assumption 
The quasi concavity of the absolute value of the utility functions {\bf (H5)} hold if $ U(z)=\bar U(z)=\log(z)$ or $ U(z)=\bar U(z)=\frac{z^\eta}{\eta}$, $\eta\in (0,1)$. 
\end{Remark}
\begin{Lemma}\label{lemme1}
We assume that the set of controls $\Ac$ is non-empty {\bf(H3)} and $x\geq v(0,0)$. 
The set $\Ac(x) $ is closed for almost everywhere convergence topology. 
%the topology of convergence in measure.
\end{Lemma}
\noindent{\bf Proof.}
We consider a sequence $(c^n,\xi^n)_n\in \Ac$ such that   
\begin{eqnarray*}
\xi^n \longrightarrow \hat \xi\,\, dP \mbox{ a.s. and }c^n_t \longrightarrow \hat c_t \,\, dt\otimes dP\mbox{ a.e. }
\end{eqnarray*}
By Fatou's lemma and using the uniform integrability of the family $\Big(\exp{(\gamma |\bar U(\xi^n)|)}\Big)_n$ ,we have
\begin{eqnarray}\label{conduni1}
E_P[ \exp{(\gamma |\bar U( \hat \xi)|)}] &\leq & \Sup_n E_P[ \exp{(\gamma |\bar U ( \xi^n)|)}]<\infty.
\end{eqnarray} 
By the uniform integrability of the family $\Big(\exp{(\gamma |\bar U( \xi^n)|)}\Big)_n$ and for a fixed $\epsilon >0$, 
there exists $\zeta >0$ such that, if $P(A)\leq \zeta$, then $\int_A \exp{(\gamma |\bar U( \xi^n)|)}dP \leq \epsilon$ which implies
\begin{eqnarray}\label{conduni2}
\int_A \exp{(\gamma |\bar U( \hat \xi)|)}dP \leq  \Lim_n\Inf \int_A \exp{(\gamma |\bar U (\xi^n)|)}dP \leq \epsilon
\end{eqnarray}
From \reff{conduni1}-\reff{conduni2}, we deduce the boundedness in $L^1(P)$ and the equi-integrability of $\hat \xi$. 
This shows $\hat \xi\in \Hc^{ad}$. Similarly $\hat c\in \Cc^{ad}$ and so $(\hat c,\hat \xi)\in \Ac$. 
By Fatou's lemma, we have
\begin{eqnarray*}
E_{P^{\nu}} \left[\hat \xi+\int_0^T\hat c_tdt
-A_T({\nu}) \right]
&\leq &
\Liminf_{n\longrightarrow \infty}E_{P^{\nu}} \left[\xi^n+\int_0^Tc^n_{t}dt
-A_T({\nu}) \right]\\
&\leq& x,
\end{eqnarray*}
which implies that $v(c,\xi)\leq x$. From the characterization \reff{vhx}, we deduce that 
$(c,\xi)\in \Ac(x)$ and so the closeness of the set $ \Ac(x)$ is proved.
\ep\\
\begin{Lemma}\label{lemme0}
We assume that the set of controls $\Ac$ is non-empty {\bf(H3)}, $x\geq v(0,0)$ and 
the quasi concavity of the absolute value of the utility functions {\bf (H5)} holds, then the set $\Ac(x) $ is convex.
\end{Lemma}
\noindent{\bf Proof.}
We take  $(c_1,\xi_1) \in \Ac(x)$, $(c_2,\xi_2) \in \Ac(x)$ and $\lambda\in [0,1]$.
From the quasi concavity of the absolute value of the utility functions {\bf (H5)}, and using Cauchy Schwartz inequality, we have for all $\gamma>0$
\begin{eqnarray*}
E_P[ \exp{(\gamma |\bar U(\lambda \xi_1 + (1-\lambda)\xi_2)|)}]
&\leq &\sqrt{ E_P[\exp{(2 \gamma |\bar U (\xi_1)|)}]}\sqrt{ E_P[\exp{(2\gamma |\bar U(\xi_2)|)}]}\\
&\leq&\Sup_{\Hc^{ad}}E_P[\exp{(2 \gamma |\bar U (\xi)|)}]<\infty.
\end{eqnarray*}
The equi-integrability of $\exp{(\gamma \bar U(\lambda \xi_1 + (1-\lambda)\xi_2))}$ holds as in the lemma \ref{lemme1} and so $\lambda \xi_1 + (1-\lambda)\xi_2 \in \Hc^{ad}$.
Similarly $\lambda c_1 + (1-\lambda)c_2 \in \Cc^{ad}$.
%\begin{eqnarray*}
%& &E^P[\int_0^T \exp{(\gamma | U(\lambda (c_{1s} + (1-\lambda)c_{2s})|)}ds]\\
%&\leq & \sqrt{ E^p[\int_0^T\exp{(2 \gamma | U (c_{1s})|)}ds]}\sqrt{ \int_0^TE^P[\exp{2\gamma |\bar U(c_{2s})|}ds]}\\
%&\leq & \Sup_{\Hc_1}  E^p[\int_0^T\exp{(2 \gamma | U (c_{s})|)}ds]<\infty,
%\end{eqnarray*}
%From Assumption {\bf (H5)(ii)}, we have
%\begin{eqnarray*}
%E_P[ \exp{(\gamma |\bar U^{'}(\lambda (\xi_1 + (1-\lambda)\xi_2)|)}]
%\leq \sqrt{ E_p[\exp{(2 \gamma |\bar U^{'} (\xi_1)|)}]}\sqrt{ E[\exp{(2\gamma |\bar U^{'}(\xi_2)|)}]}.
%\end{eqnarray*}
%Similarly, we have 
%\begin{eqnarray}
%& &E_P[\int_0^T \exp{(\gamma | U^{'}(\lambda (c_{1s} + (1-\lambda)c_{2s})|)}ds]\\
%&\leq & \sqrt{ \int_0^TE_p[\exp{(2 \gamma | U^{'} (c_{1s})|)}ds]}\sqrt{ \int_0^TE[\exp{2\gamma |\bar U^{'}(c_{2s})|}ds]}.
%\end{eqnarray}
This shows the convexity of $\Ac$. The convexity of $\Ac(x)$ follows from the convexity of $\Ac$ and Proposition \ref{dualdomin}.  
\ep\\
The following lemma shows the $L^p$ integrability of $Z^*$ for all $p\geq 1$.
\begin{Lemma}\label{Z*lp} We assume that the set of controls $\Ac$ is non-empty ${\bf(H3)}$ and $x\geq v(0,0)$. 
The density $Z^*$  of the probability measure $Q^*$ is in $L^p(P)$ for all $p\geq 1$.
\end{Lemma}
\noindent{\bf Proof.}
From the dynamics of $Y^{x,c,\xi}$ given by the equation \reff{bsdech}, we have
\begin{eqnarray*}
-\frac{1}{\beta}(Y_t^{x,c,\xi}-Y_0^{x,c,\xi})
+\frac{1}{\beta}\int_0^t(\delta_sY_s^{x,c,\xi} -\alpha U(c_s))ds
=-\frac{1}{2\beta^2}\int_0^t |Z_s^{x,c,\xi}|^2ds-\frac{1}{\beta}\int_0^tZ_s^{ x,c,\xi \prime}dW_s,
\end{eqnarray*}
and so for all $p\geq 1$, we obtain
\begin{eqnarray*}
E_P\Big[\exp\Big( p\Big(
-\frac{1}{\beta}(Y_t^{x,c,\xi}-Y_0^{x,c,\xi})
+\frac{1}{\beta}\int_0^t(\delta_sY_s^{x,c,\xi} -\alpha U(c_s))ds
\Big)\Big)\Big]=E_P\Big[(Z^*_t)^p \Big].
\end{eqnarray*}
Since $Y\in D_0^{exp}$ and $(U(c_t))_{0\leq t\leq T}\in D_1^{\exp}$, the result follows.
\ep
\\
The next result is related to the concavity and the upper semicontinuity of  the functional $(c,\xi)\longrightarrow Y_0^{x,c,\xi}$.
\begin{Proposition}\label{regularity}
We assume that the set of controls $\Ac$ is non-empty {\bf(H3)} and $x\geq v(0,0)$. 
Under the standard assumptions on the utility functions {\bf (H4)}, 
the functional $(c,\xi)\longrightarrow Y_0^{x,c,\xi}$ is strictly concave and upper semicontinuous.
\end{Proposition}
{\bf Proof.}
We fix $\lambda \in (0,1)$, $(c^1,\xi^1)\in \Ac(x)$ and $(c^2,\xi^2)\in \Ac(x)$, 
such that $P(\xi^1\neq \xi^2)>0$ or $c^1_t\neq c^2_t$ over a non null set with respect to the measure $dt\otimes dP$. 
Then by convexity of the set $\Ac(x)$, we have $(\lambda c^1+ (1-\lambda)c^2,\lambda \xi^1+ (1-\lambda)\xi^2)\in \Ac(x)$
and \\
$(Y^{x,\lambda c^1+ (1-\lambda)c^2,\lambda \xi^1+ (1-\lambda)\xi^2},Z^{x,\lambda c^1+ (1-\lambda)c^2,\lambda \xi^1+ (1-\lambda)\xi^2})$ is solution of the BSDE 
\reff{bsdech}-\reff{bsdechct} associated with $(U(\lambda c^1+ (1-\lambda)c^2),\bar U(\lambda \xi^1+ (1-\lambda)\xi^2))$. 
We set $\bar c_t=U^{-1} ( \lambda U(c_t^1)+ (1-\lambda)U(c_t^2))$, $t\in[0,T]$ and  
$\bar \xi=\bar U^{-1}(\lambda \bar U(\xi^1)+ (1-\lambda)\bar U(\xi^2))$. 
Thanks to the standard assumptions on the utility functions {\bf(H4)}, $\bar c$ and $\bar \xi$ are well-defined.
From the concavity of $U$ and $\bar U$, we have 
$U(\lambda c_t^1+ (1-\lambda)c_t^2)\geq \lambda U(c_t^1)+ (1-\lambda)U(c_t^2)=U(\bar c_t) $ $dt\otimes dP$ a.e., $t\in[0,T]$ and 
$\bar U(\lambda \xi^1+ (1-\lambda)\xi^2)\geq \lambda \bar U(\xi^1)+ (1-\lambda)\bar U(\xi^2)=\bar U(\bar \xi)$ $dP$ a.s.
The comparison theorem (See Theorem \ref{ctheorem}) yields
\begin{eqnarray}\label{concav1}
Y_t^{x,\lambda c^1+ (1-\lambda)c^2,\lambda \xi^1+ (1-\lambda)\xi^2}\geq Y_t^{x,\bar c,\bar \xi}\,\,dt\otimes dP\,\,\mbox{ a.e.},\, t\in[0,T].
\end{eqnarray}
From the definition of $Y_t^{x,\bar c,\bar \xi}$ (See equation \reff{robust})
\begin{eqnarray}\label{concav2}
Y_t^{x,\bar c,\bar \xi}&=& \mbox{ess}\Inf_{Q\in \Qc_f}\Big(
\frac{1}{S_t^\delta}
E_Q\Big[ \int_{t}^{T}\alpha S_s^\delta U(\bar c_s)ds+ \bar \alpha S_T^\delta \bar U(\bar \xi) \Big |\Fc_t \Big]
+\beta E_Q\Big[\Rc_{t,T}^{\delta}(Q)
|\Fc_t \Big]
\Big)\nonumber\\
&\geq&\lambda \mbox{ess}\Inf_{Q\in \Qc_f}\Big(
\frac{1}{S_t^\delta}
E_Q\Big[ \int_{t}^{T}\alpha S_s^\delta U( c^1_s)ds+ \bar \alpha S_T^\delta \bar U( \xi^1) \Big |\Fc_t \Big]
+\beta E_Q\Big[\Rc_{t,T}^{\delta}(Q)
|\Fc_t \Big]
\Big)\nonumber\\
&+& (1-\lambda)\mbox{ess}\Inf_{Q\in \Qc_f}\Big(
\frac{1}{S_t^\delta}
E_Q\Big[ \int_{t}^{T}\alpha S_s^\delta U( c^2_s)ds+ \bar \alpha S_T^\delta \bar U( \xi^2) \Big |\Fc_t \Big]
+\beta E_Q\Big[\Rc_{t,T}^{\delta}(Q)
|\Fc_t \Big]
\Big)\nonumber\\
&=& \lambda Y_t^{x, c^1,\xi^1} +(1-\lambda) Y_t^{x,c^2,\xi^2}.
\end{eqnarray}
From \reff{concav1} and \reff{concav2}, we deduce that
\begin{eqnarray}
Y_0^{x,\lambda c^1+ (1-\lambda)c^2,\lambda \xi^1+ (1-\lambda)\xi^2} \geq  \lambda Y_0^{x, c^1,\xi^1} +(1-\lambda) Y_0^{x,c^2,\xi^2}.
\end{eqnarray}
From the strict concavity of the utility functions, we have 
$P(\bar U(\lambda \xi^1+(1-\lambda) \xi^2)>\lambda \bar U(\xi^1)+(1-\lambda) \bar U(\xi^2))>0$ 
or $U(\lambda c^1_t+(1-\lambda) c^2_t)>\lambda U( c^1_t)+(1-\lambda) U(c^2_t)$ over a non null set with respect to the measure $dt\otimes dP$.
The comparison theorem (See Theorem \ref{ctheorem}) yields
\begin{eqnarray*}
Y_0^{x,\lambda c^1+ (1-\lambda)c^2,\lambda \xi^1+ (1-\lambda)\xi^2} > \lambda Y_0^{x, c^1,\xi^1} +(1-\lambda) Y_0^{x,c^2,\xi^2}.
\end{eqnarray*}
This shows the strict concavity of $(c,\xi)\longrightarrow Y_0^{x,c,\xi}$.\\
We turn to the upper semicontinuity of $Y_0^{x,c,\xi}$. 
Let $(c^n,\xi^n)\in \Ac(x)$ such that $c^n_t\longrightarrow c_t$, $dt\otimes dP$ a.e., $t\in[0,T]$ and $\xi^n\longrightarrow \xi$ dP a.s. 
From Lemma \ref{lemme1}, we have $(c,\xi)\in \Ac(x)$. We set $\tilde c^n=\Sup_{m\geq n}c^m$ and $\tilde \xi^n=\Sup_{m\geq n}\xi^m$. 
Then $\tilde \xi^n\searrow \xi$ $dP$ a.s. and $\tilde c^n_{t}\searrow c_t $,
$0\leq t\leq T$, $dt\otimes dP$ a.e. when $n$ goes to infinity. From Proposition \ref{monotonieYn}, we deduce that 
then $Y_0^{x, \tilde c^n,\tilde \xi^n}\searrow Y_0^{x, c,\xi}$. 
On the other hand, we have $\tilde \xi^n\geq \xi^n $ and $\tilde c^n\geq c^n$ and so the comparison theorem (See Theorem \ref{ctheorem}) yields 
$Y_0^{x, \tilde c^n,\tilde \xi^n} \geq  Y_0^{x, c^n,\xi^n}$ which implies $\Lim_n Y_0^{x, \tilde c^n,\tilde \xi^n} \geq  \Lim_n \Sup Y_0^{x, c^n,\xi^n}$. 
This shows that $Y_0^{x, c,\xi} \geq  \Lim_n \Sup Y_0^{x, c^n,\xi^n}$ and so $(c,\xi)\longrightarrow Y_0^{x,c,\xi}$ upper semicontinuous.
\ep \\
The next lemma shows the boundedness of the value function.
\begin{Lemma}\label{lemuni} We assume that the discounting factor is bounded {\bf (H1)}, the set of controls $\Ac$ is not reduced to the null strategy {\bf(H3)} and $x\geq v(0,0)$. we have
\begin{eqnarray}\label{finitude}
\Sup_{(c,\xi)\in \Ac(x)}Y_0^{x,c,\xi}<\infty.
\end{eqnarray}
\end{Lemma}
{\bf Proof.}
From the definition of $Y_0^{x,c,\xi}$ and using the boundedness on the discounting factor {\bf (H1)}, we have
\begin{eqnarray*}
Y_0^{x,c,\xi}&\leq & E_{P}[\int_0^T \alpha S^{\delta}_s U(c_s)ds+ \bar \alpha S^{\delta}_T \bar U(\xi) ]\\
%&+& \beta  E_{P}[\int_{0}^{T}\delta_s S_s^\delta  Z_s^{*} \log Z_s^{*} ds+
%S_T^\delta Z_T^{*}\log  ]\\
&\leq & C \Big(E_{P}[\int_0^T |U(c_s)|ds+|\bar U(\xi)| ]\Big).
%+  E_{P}[\int_{0}^{T} |Z_s^{*}\log Z_s^{*}| ds+ |Z_T^{*}\log Z_T^{*}| ],
\end{eqnarray*}
Since for all $y\geq 0$, we have $|\bar U(y)|\leq C(1+ \exp{|\bar U (y)|})$ and $|U(y)|\leq C(1+\exp{|U (y)|})$, then the result follows from the uniform integrability
of the families \reff{unifinteg1} and \reff{unifinteg3}.
\ep\\\\
Our next result is the existence of a unique solution to the problem \reff{reward}.
The uniqueness follows since $(c,\xi)\longrightarrow Y_0^{x,c,\xi}$ is strictly concave.
\begin{Theorem} \label{main1} 
We assume that the discounting factor is bounded {\bf (H1)}, the set of controls $\Ac$ is non-empty {\bf(H3)},
the utility functions satisfy the usual conditions{\bf (H4)}, 
the absolute value of the utility functions is quasi-concave  {\bf (H5)}, 
and $x\geq v(0,0)$.
There exists a unique solution
$(c^*,\xi^*)\in  \Ac(x)  $ of \reff{reward}.
 \end{Theorem}
{\bf Proof.}
Let $(c^n,\xi^n)_{n\in \N}\in  \Ac(x)$ be a maximizing sequence of the problem \reff{reward} i.e.
\begin{eqnarray}\label{solution}
\Lim_{n \to \infty} Y_0^{x,c^n,\xi^n}  \; = \; \Sup_{(c,\xi)\in \Ac(x) } Y_0^{x,c,\xi},
\end{eqnarray}
which is finite by Lemma \ref{lemuni}.\\
Since $\xi^n \geq 0\;\; dP\,a.s$ and $c^{n}_{t} \geq 0\;\; dt \otimes dP\,a.s$,
then by Lemma A.1.1 of Delbaen and Schachermeyer \cite{delsch94}, there exists a sequence
$(\hat c^n, \hat \xi^n) \in \mbox{conv }((c^n,\xi^n),(c^{n+1},\xi^{n+1}),...)$ such that  $(\hat c^n, \hat \xi^n)$
converges almost surely to $(c^*, \xi^*)\in  \tilde \Cc \times L^0_+(\Fc_T)$. By Lemmas \ref{lemme1}-\ref{lemme0}, 
we have $(\hat c^n, \hat \xi^n)\in \Ac(x)$ and $(c^*, \xi^*)\in \Ac(x) $. From Proposition \ref{regularity},the functional
$(c,\xi)\longrightarrow Y^{x,c,\xi}$ is concave and so
\begin{eqnarray*}
\Sup_{(c,\xi)\in  \Ac(x)}Y_0^{x,c,\xi}
&\leq& \Limsup_{n \to \infty}Y_0^{x,\hat c^n,\hat \xi^n}.
\end{eqnarray*}
From Proposition \ref{regularity},
the functional $(c,\xi)\longrightarrow Y^{x,c,\xi}_0$ is upper semicontinuous and so 
\begin{eqnarray*}
\Limsup_{n \to \infty}Y_0^{x,\hat c^n,\hat \xi^n}\leq Y_0^{x,c^*,\xi^*}.
\end{eqnarray*}
Therefore $(c^*,\xi^*)$ solves \reff{reward}.
\ep
%\begin{Remark} El Karoui et al. \cite{elk} use another approach  based on the functional analysis
%to prove the existence of an optimal
%solution of the  optimization problem \reff{reward2}.
%They penalise their criterion
%in order to obtain the strong convexity property. They prove  the existence of an
%optimal strategy for the penalized problem. Then, they show that the sequence
%of the optimal strategy in uniformly bounded for the weak topology and so the
%existence of an optimal strategy of the problem \reff{reward2} follows.
%\end{Remark}
\section{Duality and Dynamic Maximum Principle}
The aim of this section is to provide a description of the solution structure to problem \reff{reward} via the dual formulation. 
%The heart of the dual approach in the classical setting, when the creterion is taken under the historical 
%probability measure, is to find a saddle point for the Lagrangian and apply a minimax theorem in the infinite dimensional case.
%It is natural to use the conjugate function of U and $\bar U$. 
In fact, to solve an investment problem when the underlying model is known, and by using the definition of the conjugate function of $\bar U$ denoted by $\tilde U$, we have
\begin{eqnarray*}
\bar U(X_T^{x,H})\leq \tilde U(yZ_T^{\nu})+yZ_T^{\nu}X_T^{x,H}.
\end{eqnarray*}
If $(Z_t^{\nu}X_t^{x,H})_{t\in [0,T]}$ is a supermartingale under $P$, we have $E[Z_T^{\nu}X_T^{x,H}]\leq x$, which implies
\begin{eqnarray*}
\Sup_{H}E[\bar U(X_T^{x,H})]\leq \Inf_{\nu}E[\tilde U(yZ_T^{\nu})]+xy.
\end{eqnarray*}
If we find $H^*$ and $\nu^*$ such that we have equality in the latter equation for some $y^*$, then $\nu^*$ is the solution of the dual problem.
In our case, the criterion in taken 
under $Q^*$ and the use of the conjugate functions is not appropriate. In fact we have 
\begin{eqnarray*}
Z_T^*\bar U(X_T^{x,H})\leq Z^*_T\tilde U(yZ_T^{\nu})+yZ^*_TZ_T^{\nu}X_T^{x,H},
\end{eqnarray*}  
and so the supermartingale property  of $(Z_t^*Z_t^{\nu}X_t^{x,H})_{t\in [0,T]}$ does not hold in general.
We will use the arguments of the duality differently.
%First, we will show 
%the existence of an optimal probability measure solution of the dual problem.
%that the budget constraint is satified with equality.
%To prove such result, we need convex duality arguments. 
First, we will show that there exists a probability measure $\tilde P^{^*}$ equivalent to the probability measure $P$ solution of the problem
\begin{eqnarray}\label{du}
v(c^*,\xi^*)=\sup_{P^\nu\in\Pc^0} 
E_{P^\nu} \left[\xi^*+\int_0^Tc^*_tdt
-A_T(\nu) \right].
\end{eqnarray}
Then, we will show that the budget constraint is satisfied with 
equality which is a consequence from the strict concavity of the utility functions and the comparison theorem.
We start with the following lemma.
\begin{Lemma}\label{conv}
The set of probability measures $\Pc^0$ is convex and the function $P^{\nu}\longrightarrow E_{P^{\nu}}[A_T(\nu)]$ is convex.
\end{Lemma}
{\bf Proof.}\\
Let $P^{\nu^1},\,\,P^{\nu^2}\,\,\in \Pc^{0},\,\, Z^{\nu^1},\,\,Z^{\nu^2}$ their density processes, $\alpha\, \in\,[0,1]$
and denote by $P^{\tilde \nu }\sim P $ the probability measure $P^{ \tilde \nu }=\alpha P^{\nu^1}+(1-\alpha)P^{\nu^2}$ and by $Z^{\tilde \nu }$
its density process. Consider the process $A^{P^{\tilde \nu }}$ defined by
\begin{eqnarray}\label{deup}
A^{P^{\tilde \nu }}_t=\alpha \int_0^t\frac{Z^{\nu^1}_u}{Z_u^{\tilde \nu }}dA_u(\nu^1)+(1-\alpha)\int_0^t\frac{Z^{\nu^2}_u}{Z_u^{\tilde \nu }}dA_u(\nu^2)
\,\,\,\,\,\,0\leq t\leq T.
\end{eqnarray}
From \reff{propupper}, we have $ A_T(\tilde \nu )\leq A_T^{P^{\tilde \nu}}$, which implies that, 
\begin{eqnarray*}
 A_T(\tilde \nu )\leq \alpha \int_0^t\frac{Z^{\nu^1}_u}{Z_u^{\tilde \nu }}dA_u(\nu^1)+(1-\alpha)\int_0^t\frac{Z^{\nu^2}_u}{Z_u^{\tilde \nu }}dA_u(\nu^2)
\,\,\,\,\,\,0\leq t\leq T.
\end{eqnarray*}
Since, $0\leq \frac{Z^{\nu^i}_u}{Z_u^{\tilde \nu }}\leq 1$, for $i=1,2,$ we obtain,
 \begin{eqnarray*}
 A_T(\tilde \nu )&\leq& \alpha \int_0^t dA_u(\nu^1)+(1-\alpha)\int_0^t dA_u(\nu^2)
\,\,\,\,\,\,0\leq t\leq T.\\
&\leq&  \alpha A_T(\nu^1)+(1-\alpha) A_T(\nu^2).
\end{eqnarray*}
We have $\nu^i \in \Nc(\tilde K)$, which implies that $A_T(\nu^i)<+\infty$ for  $i=1,2,$ and so  
$ A_T(\tilde \nu )< +\infty$.\\
From the convexity of the functions $z\longrightarrow z^{\eta}$ where $\eta> 1$ and using the definition of $\Pc^0$, we have 
\begin{eqnarray*}
E[(Z^{\tilde \nu}_T)^{\eta}]&\leq& \alpha E[(Z^{\nu^1}_T)^{\eta}]+ (1-\alpha)E[(Z^{\nu^2}_T)^{\eta}]\\
&\leq & \Sup_{\nu}E[(Z^{\nu}_T)^{\eta}] <\infty.
\end{eqnarray*}
Similarly, from the convexity of the function $z\longrightarrow z^{1-\bar \eta}$, where $\bar\eta> 1$, we have 
$E[(Z^{\tilde \nu}_T)^{1-\bar \eta}] \leq  \Sup_{\nu}E[(Z^{\nu}_T)^{1-\bar \eta}] <\infty$.
To check the equi-integrability point, for all $I=[t_0,t_1] \subset [0,T]$, we have
\begin{eqnarray*}
\int_{t_0}^{t_1}\delta^{supp} (\tilde \nu_t)dt &=&A_{t_1}(\tilde \nu )-A_{t_0}(\tilde \nu )\leq A_{t_1}^{P^{\tilde \nu}}-A_{t_0}^{P^{\tilde \nu}}\\
&\leq & \int_{t_0}^{t_1}\delta^{supp} ( \nu^1_t)dt+\int_{t_0}^{t_1}\delta^{supp} ( \nu^2_t)dt.
\end{eqnarray*}  
The first inequality is deduced from \reff{propupper}  and 
the second one is deduced from equation \reff{deup} and by using the equality $Z_T^{ \tilde \nu }=\alpha Z_T^{\nu 1}+(1-\alpha)Z_T^{\nu 2}$. \\
Let $\displaystyle{\epsilon > 0,\,\,\mbox{there exists}\,\, \epsilon_1>0\,\,\mbox{such that if}\,\,\lambda(I) \leq \epsilon_1,\,\, \mbox{then, we have}\,\,\int_I \delta^{supp} ( \nu^1_t)dt \leq \frac{\epsilon}{2}},$
$\displaystyle{\mbox{and}}$\\
$\displaystyle{\mbox{ there exists}\,\, \epsilon_2>0\,\,\mbox{such that if}\,\,\lambda(I)\leq \epsilon_2,\,\, \mbox{then, we have}\,\,\int_I \delta^{supp} ( \nu^2_t)dt \leq \frac{\epsilon}{2}}.$\\
We take $\displaystyle{\epsilon =\inf(\epsilon_1, \epsilon_2)\,\, \mbox{, it's follows that }\,\,\int_I \delta^{supp} (\tilde \nu_t)dt \leq\epsilon.}$\\
This shows the equi-integrability property of $\delta^{supp} (\tilde \nu .)$ with respect to the Lebesgue measure and so the convexity of  $\Pc^{0}$. \\
From the inequality $ A_T(\tilde \nu )\leq A_T^{P^{\tilde \nu}}$, we have
\begin{eqnarray*}
\E_{P^{\tilde \nu}}\big[A_T(\tilde \nu) \big]\leq  \alpha \E_{P^{\nu 1}} \big[ A_T(\nu 1) \big]+ (1-\alpha) \E_{P^{\nu 2}} \big[ A_T(\nu_2) \big],
\end{eqnarray*}
and so we deduce the convexity of the function $P^{\nu }\in \Pc^{0}\longrightarrow \E^{P^{\nu}}\big[A_T(P^{\nu}) \big]$.
\ep\\

The following Theorem shows the existence of a probability measure $\tilde{P}^{*}$ equivalent to the probability measure $P$ solution of the problem
\reff{du}. The density process of $\tilde{P}^{*}$ with respect to $P$ is the $P$-martingale
$\tilde{Z}^*=(\tilde{Z}^*_t)_{0\leq t\leq T}$ with
\begin{eqnarray}\label{Ztilde}
\tilde{Z}_t^*=E_P\Big[ \frac{d\tilde{P}^{*}}{dP}\Big |\Fc_t\Big]; t\in [0,T], \,\,dt\otimes dP,a.e.
\end{eqnarray}
%Here, we assume that the supremum over all probabilty measures in $\Pc^0$ of the budget 
%constraint is attained under $\tilde P^*$. The criterion is not concave in $\nu$ and the set of $\nu$ such that $P^{\nu}\in \Pc^0$ is not convex. 
%The coercivity property is not satisfied.
% Also, it is not clear how we extend appropriately the set $\Pc^0$ to 
%prove the exsitence of $\tilde P^*$.
We shall assume the translation stability on the set of admissible strategies i.e. \\ 
{\bf (H6)} If $(c,\xi)\in \Ac(x)$, then $(c+ \alpha^{1},\xi +\alpha^{2})\in \Ac(x)$ for any $\alpha^{1}>0$ and $\alpha^{2}>0$.
\begin{Remark}
The translation stability on the set of admissible strategies {\bf (H6)} is satisfied if the utility functions are subadditive.
\end{Remark}
\begin{Theorem}\label{bd}
We fix $x\geq v(0,0)$. We assume that the discounting factor is bounded {\bf (H1)}, the set of controls $\Ac$ is non-empty {\bf(H3)}, 
the utility functions satisfy the usual conditions {\bf (H4)}, the absolute value of the utility functions is quasi-concave {\bf (H5)},  and the translation stability on the
set of admissible strategies {\bf (H6)} holds. Then, 
there exists a probability measure $\tilde P^*\in \Pc^0 $  such that 
\begin{eqnarray}\label{budgetc1} 
\sup_{P^\nu\in\Pc^{0}} 
E_{P^\nu} \left[\xi^*+\int_0^Tc^*_tdt
-A_T(\nu) \right]&=&E_{\tilde P^*} \left[\xi^*+\int_0^T c^*_tdt
-A_T( \nu^*) \right],
\end{eqnarray}
and the budget constraint is satisfied with equality i.e.
\begin{eqnarray}\label{budgetc} 
\sup_{P^\nu\in\Pc^{0}} 
E_{P^\nu} \left[\xi^*+\int_0^Tc^*_tdt
-A_T(\nu) \right]=x.
\end{eqnarray}
\end{Theorem}
{\bf Proof.}\\
Let $F(P^{\nu})$ and $G(P^{\nu})$ defined by the following functionals:
\begin{eqnarray*}
F(P^{\nu})=\xi^*+\Int_0^T c_t^*dt-A_T(\nu), 
\end{eqnarray*}
and
\begin{eqnarray*}
G(P^{\nu})=E_{P^{\nu}}\big[F(P^{\nu}) \big].
\end{eqnarray*}
%The main result of this theorem is that of maximizing $G(P^{\nu})$ over $P^{\nu} \in \Pc^{0}$ has a %solution $\tilde{P}^{*}$ is even equivalent to $P$.\\
$\star$ \underline{First step}: Let $(P^{\nu_n})_{n \in \N}$ be a sequence in $\Pc^{0}$ such that:
\begin{eqnarray*}
\lim_{n\longrightarrow +\infty} G(P^{\nu_n})=\sup_{P^{\nu} \in \Pc^{0}}G(P^{\nu })\,\,<\,\infty.
\end{eqnarray*}
and denote by $Z^n=Z^{P^{\nu_n}}$ the corresponding density process.
Since each $Z^n_T \geq 0$, it's follows from Komlos' theorem that there exists a sequence $(\bar{Z}^n_T)_{n\in \N}$
with $(\bar{Z}^n_T)_{n\in \N}\in \mbox{conv}(Z^n_T,Z^{n+1}_T,\dots)$ for each $n \in \N$ and such that $(\bar{Z}^n_T)$
converge $P.a.s$ to some random variable $(\bar{Z}^{\infty}_T)$ , which is then also non-negative but may take value $+\infty$.\\
Because $\Pc^{0}$ is  convex, each $\bar{Z}^n_T$ is again associated to some $\bar{P}^{n}$ which is in $\Pc^{0}$. 
By de la Vall\'ee-Poussin's criterion, $(\bar{Z}^n_T)_{n\in \N}$ is uniformly integrable and therefore converges in $L^1(P)$. This implies that 
$\Lim_{n\longrightarrow \infty}E_P[\bar{Z}^n_T]=E_P[\bar{Z}^{\infty}_T]=1$ and so $d\bar{P}^{\infty}=\bar{Z}^{\infty}_T\,\, \mbox{dP}\,\,$ defines a probability
measure which is absolutely continuous with respect to $P$.\\
We define the following stopping time: $\tau^{\infty}=\inf \{t\geq 0\mbox{ s.t. } \bar Z_t^{\infty}=0\}$. From the continuity property of the process $\bar Z^{\infty}$,
on the set $A:=\{\tau^{\infty}\leq T\}$, we have $\bar Z^{\infty}_{\tau^{\infty}}=0$ dP a.s. Using the martingale property of $\bar{Z}^{\infty}$, we deduce that   
$\bar{P}^{\infty}(A)=E_P[\bar Z^{\infty}_{\tau^{\infty}}1_A]=0$. From the inequality
\begin{eqnarray*}
|\bar P^n(A)-\bar P^{\infty}(A)|\leq E_P[|\bar Z^n_T-\bar Z^{\infty}_T|1_A] \leq E_P[|\bar Z^n_T-\bar Z^{\infty}_T|],
\end{eqnarray*} 
and since $\bar Z^n_T$ converges to $\bar Z^{\infty}_T$ in $L^1(P)$, we deduce that $\Lim_{n\longrightarrow \infty}\bar P^n(A)=0$. 
Since $P$ and $\bar P^n$ are equivalent probability measures, we have
\begin{eqnarray*}
P(A) = E_{\bar P^n}[\frac{1}{\bar Z^n_T}1_A]&\leq& (E_{ \bar P^n}[(\bar{Z}^n_T)^{-\bar \eta}])^\frac{1}{\bar \eta}\bar P^n(A)^{1-\frac{1}{\bar \eta}}\\
&\leq & (E_{ P}[(\bar{Z}^n_T)^{1-\bar \eta}])^\frac{1}{\bar \eta}\bar P^n(A)^{1-\frac{1}{\bar \eta}}.
\end{eqnarray*} 
From the definition of the set $\Pc^0$, and since $E_{ P}[(\bar{Z}^n_T)^{1-\bar \eta}]$ is finite, there exists a positive constant C such that 
$P(A)\leq C \bar P^n(A)^{1-\frac{1}{\bar \eta}}$. Sending $n$ to infinity, we conclude that $P(A)=0$ and so  $\bar P^{\infty}$ is a probability measure which is equivalent to $P$.\\
$\star$ \underline{Second step}: We will show that $G(\bar{P}^{\infty})\geq G(P^{\nu}) $ for all  $P^{\nu}\in \Pc^{0}$.\\ %and $\bar{P}^{\infty}\in \bar{\Pc}^{0}$.\\
%on $\Pc^{0}$ \\
Since we know that $ (\bar{Z}^n_T)_n $ converges to $ \bar{Z}^{\infty} $ in $ L^1(P),$ the Doob's maximal inequality
$$P[\sup\limits_{0\leq t \leq T} \mid \bar{Z}^{\infty}_t-\bar{Z}^n_t \mid\geq \epsilon]\leq \frac{1}{\epsilon} E_P[\mid \bar{Z}^{\infty}_T-\bar{Z}^n_T \mid]$$
implies that  $(\sup\limits_{0\leq t \leq T} \mid \bar{Z}^{\infty}_t-\bar{Z}^n_t \mid)_{n\in\N}$
converges to $0$ in $P$-probability.\\
Going to a sub-sequence, still denoted by $(\bar{Z}^n)_{n\in \N}$, we can assume that \\
$(\sup\limits_{0\leq t \leq T} \mid \bar{Z}^{\infty}_t-\bar{Z}^n_t \mid)_{n\in\N}$ converges to $0$  $P$-a.s. 
\\Let $M_t^n:=\sup\limits_{0\leq s \leq t} \mid \bar{Z}^{\infty}_s-\bar{Z}^n_s \mid $ and $(\tau_n)$ a sequence of stopping time defined by
\begin{equation*}
\tau_n=\left\{\begin{array}{ccc}
         \inf\{t\in [0,T); M^n_t\geq 1\} & \textrm{if}  & \{t\in [0,T); M^n_t \geq 1\}\neq \emptyset \\
         T &  & \textrm{otherwise}
       \end{array}\right..
\end{equation*}
Since $M_{\tau_n}^n$ is bounded by $M_T^n \wedge 1$ then $M_{\tau_n}^n$ converges almost surely to $0$ 
and, by the dominated convergence theorem, converges to $0$ in $L^1(P).$
Then, using Burkholder Davis Gundy inequality $\langle \bar{Z}^{\infty}-\bar{Z}^n \rangle_{\tau_n}^\frac{1}{2}$ 
converges to $0$ in $L^1(P)$ and a fortiori in probability.
\\ As,  $\langle \bar{Z}^{\infty}-\bar{Z}^n \rangle_T= \langle \bar{Z}^{\infty}-\bar{Z}^n \rangle_{\tau_n}\textbf{1}_{\{\tau_n=T\}}
+\langle \bar{Z}^{\infty}-\bar{Z}^n \rangle_T\textbf{1}_{\{\tau_n<T\}}$, then for all $\varepsilon>0,$
\begin{equation*}
\begin{split}
P(\langle \bar{Z}^{\infty}-\bar{Z}^n \rangle_T\geq \varepsilon) & \leq P(\langle \bar{Z}^{\infty}-\bar{Z}^n \rangle_{\tau_n}\textbf{1}_{\{\tau_n=T\}}\geq \varepsilon)+ P(\langle \bar{Z}^{\infty}-\bar{Z}^n \rangle_T\textbf{1}_{\{\tau_n<T\}}\geq \varepsilon)
\\ & \leq P(\langle \bar{Z}^{\infty}-\bar{Z}^n \rangle_{\tau_n}\geq \varepsilon)+P(\tau_n<T).
\end{split}
\end{equation*}
From the convergence in probability of $(\langle \bar{Z}^{\infty}-\bar{Z}^n \rangle_{\tau_n})_n$, 
we have $\lim\limits_{n\rightarrow +\infty}P(\langle \bar{Z}^{\infty}-\bar{Z}^n \rangle_{\tau_n}\geq \varepsilon)=0.$
Since $M^n$  is an increasing process, we have
$$P(\tau_n<T)=P(\{\exists t\in[0,T)\;\; s.t\;\; M^n_t\geq 1\})\leq P(\{ M^n_T\geq 1\}).$$
Since $M_T^n$ converges in probability to $0$, we have $P(\{ M^n_T\geq 1\}) \underset{n\rightarrow +\infty}{\longrightarrow} 0$. 
Then $\lim\limits_{n\rightarrow +\infty}P(\tau_n<T)=0,$
and consequently $\lim\limits_{n\rightarrow +\infty}P(\langle \bar{Z}^{\infty}-\bar{Z}^n \rangle_T\geq \varepsilon)=0$ i.e. 
$(\langle \bar{Z}^{\infty}-\bar{Z}^n \rangle_T)_n $ converges in probability to $0.$
We can extract a sub-sequence denoted also by $\bar{Z}^n$ such that  $(\langle \bar{Z}^{\infty}-\bar{Z}^n \rangle_T)_n $ converges  almost surely to $0.$\\
On the other hand, we have $\bar P^{\infty}$ is equivalent to $P$, 
which implies that $\bar Z_T^{\infty}>0$ $P$ $a.s.$ and we have $\bar \nu ^{\infty}\in \Nc$ such that 
$\bar Z^{\infty}_T= \Ec_T \left( - \int (\theta + \sigma^{-1} \bar \nu^{\infty})^{\prime}dW \right)$. It yields that, when $n$ goes to infinity, we have 
$$\langle \bar{Z}^{\infty}-\bar{Z}^n \rangle_T=\int_0^T(\bar{Z}_t^{n}(\theta_t+\sigma^{-1}_t{\bar \nu^n}_t)
-\bar{Z}^{\infty}_t(\theta_t+\sigma^{-1}_t{\bar \nu^{\infty}}_t))^2du \longrightarrow 0.$$
Since $\bar{Z}^n\longrightarrow \bar{Z}^{\infty} dt\otimes dP$-a.e, we have $\bar{\nu}^n$ converges to $\bar{\nu}^\infty$ $dt\otimes dP$- a.e.\\
The continuity of the support function $\delta^{supp}$ (see Assumption \ref{fonctsupp}) yields that $\delta^{supp} (\bar{\nu}^n)$ converges to $\delta^{supp}(\bar{\nu}^\infty)$ 
$dt\otimes dP$- a.e.
From the definition of the set $\Pc^0$, $\big(\delta^{supp} (\bar{\nu}^n)\big)_n$ is equi-integrable with respect to the Lebesgue measure, and so, we have
\begin{eqnarray}\label{convPS}
\int_0^T \delta^{supp}(\bar \nu_t^n)dt \longrightarrow \int_0^T\delta^{supp}(\bar \nu_t^{\infty})dt\,\,,   \text{P-a.s.}\,\,\,\mbox{when}\,\, n \rightarrow \infty.  
\end{eqnarray}
This yields that, when $n$ goes to infinity,
\begin{eqnarray}\label{cvps}
\bar{Z}^{n}_T F(\bar{P}^{n})\longrightarrow \bar{Z}^{\infty}_T F(\bar{P}^{\infty})\,\,,\text{P-a.s.}
\end{eqnarray}
From de la Vall\'ee Poussin's criterion, we deduce the uniform integrability of the family $(\bar{Z}^{n}_T F(\bar{P}^{n}))_n$. 
In fact, for a fixed $\eta^{'}$ satisfying $\eta>\eta^{'}>1$ and $\frac{\eta^{'}\eta}{\eta-\eta^{'}}=3$, the Cauchy Schwartz inequality implies
\begin{eqnarray}\label{uni}
E_P\bigg[(\bar{Z}^{n}_T F(\bar{P}^{n}))^{\eta^{'}}\bigg]
&=&E_P\bigg[(\bar{Z}^{n}_T (\xi^*+\int_0^T c^*_tdt-\int_0^T\delta^{supp}(\bar \nu_t^{n})dt))^{\eta^{'}}\bigg]\nonumber\\
&\leq & E_P\bigg[(\bar{Z}^{n}_T)^{\eta}\bigg]^{\frac{\eta^{'}}{\eta}}E_P\bigg[(\xi^*+\int_0^T c^*_tdt 
-\int_0^T\delta^{supp}(\bar \nu_t^{n})dt)^{3}\bigg]^{1-\frac{\eta^{'}}{\eta}}\nonumber\\
&\leq& E_P\bigg[(\bar{Z}^{n}_T)^{\eta}\bigg]^{\frac{\eta^{'}}{\eta}}E_P\bigg[(\xi^*+\int_0^T c^*_tdt)^{3}\bigg]^{1-\frac{\eta^{'}}{\eta}}< \infty,
\end{eqnarray}
where the second inequality is deduced from the non decreasing property of the function \\$z\longrightarrow z^3$. From \reff{cvps} and \reff{uni}, we obtain the convergence in $L^{1}(P)$ of the sequence $(\bar{Z}^{n}_T F(\bar{P}^{n}))_n$, which yields
 \begin{eqnarray*}
E_P\bigg[\bar{Z}^{\infty}_T F(\bar{P}^{\infty})\bigg]= \Lim _{n \rightarrow \infty }
E_P\bigg[\bar{Z}^{n}_T F(\bar{P}^{n}) \bigg].
\end{eqnarray*}
This shows that
\begin{eqnarray*}
 \displaystyle{
G(\bar{P}^{\infty})=E_P\bigg[\bar{Z}^{\infty}_T F(\bar{P}^{\infty})\bigg] =\Lim_{n \rightarrow \infty }
G(\bar{P}^{n})}.
\end{eqnarray*}
From the convexity of the function  $P^{\nu}\longrightarrow E_{P^{\nu}}[A_T(\nu)]$ (see Lemma \ref{conv}), we deduce the concavity of $P^{\nu}\longrightarrow G(P^{\nu})$, which implies
\begin{eqnarray*}
\Lim_{n \rightarrow \infty }
G(\bar{P}^{n}) \geq  \Sup_{P^{\nu}\in \Pc^{0}} G(P^{\nu}),
\end{eqnarray*}
%which proves that $\bar{P}^{\infty}$ is indeed optimal.\\
and so $G(\bar{P}^{\infty}) \geq  \Sup_{P^{\nu}\in \Pc^{0}} G(P^{\nu})$. We denote by $\nu^*=\bar{\nu}^{\infty}$ and $\tilde P^*$ the probability measure associated with $\tilde{Z}^{*}$, i.e. $\tilde{P}^*=\bar{P}^{\infty}$.\\\\
$\star$ \underline{Third step}: we show that the budget constraint is satisfied with equality.\\
We assume that 
\begin{eqnarray*}
\sup_{P^\nu\in\Pc^{0}} 
E_{P^\nu} \left[\xi^*+\int_0^Tc^*_tdt
-A_T(\nu) \right]= l<x.
\end{eqnarray*} 
From the characterization \reff{vhx}, we deduce that there exists $H^*\in \tilde \Hc$ such that $\xi^*\leq X_T^{l,c^*,H^*}$ where 
\begin{eqnarray*}
X_t^{l,c^*,H^*}=l+ \int_0^tH_s^*dS_s-\int_0^t c_s^*ds,\,\, dt\otimes dP\,{ a.e.}, \,t\in [0,T].
\end{eqnarray*}
We denote by $\displaystyle{\tilde c_t=c^*_t+\frac{x-l}{T}}$ $dt\otimes dP\,{ a.e.}, \,t\in [0,T]$. Then 
\begin{eqnarray*}
X_T^{x,\tilde c,H^*}&=&x+ \int_0^T H_s^*dS_s-\int_0^T \tilde c_sds\\
&=&l+ \int_0^T H_s^*dS_s-\int_0^T c_s^*ds\\
&=&X_T^{l,c^*,H^*}.
\end{eqnarray*}
%which implies $X_T^{x,\tilde c,H^*}=X_T^{l,c^*,H^*}$. 
Under the translation stability on the set of admissible strategies {\bf (H6)}, 
$(\tilde c,\xi^*)$ satisfies \reff{unifinteg1}-\reff{unifinteg3} and 
\begin{eqnarray*}
\sup_{P^\nu\in\Pc^{0}} 
E_{P^\nu} \left[\xi^*+\int_0^T\tilde c_tdt
-A_T(\nu) \right]\leq x,
\end{eqnarray*} 
which implies $(\tilde c,\xi^*)\in \Ac(x)$ (See characterization \reff{vhx}). From the comparison theorem (See Theorem \ref{ctheorem} ), we have
$Y_0^{x,\tilde c,\xi^*}\geq Y_0^{x,c^*,\xi^*} =\Sup_{(c,\xi)\in \Ac(x)}Y_0^{x,c,\xi}$ which contradicts the unicity of the optimal 
strategy $(c^*,\xi^*)$ and so $l=x$, the equality \reff{budgetc} holds. \\
$\star$ \underline {Fourth step}: We prove that $\tilde {P}^{*}$ is indeed optimal. 
Since $ \tilde{P}^*$ is equivalent to $P$, it is clear that, $\nu^*\in \Nc$. From Assumption (\ref{fonctsupp}) $\tilde{K}$ is a closed set, and so $\bar{\nu}_t^n \longrightarrow \nu^*_t $, $dt \otimes  dP\,\,a.e.$, implies that  $\nu_t^*\in \tilde{K}\,\, dt \otimes  dP $.\\
We have $E_{\tilde{P}^*}[\xi^* +\int_0^T c^*_s ds - A_T( \nu^*)]\geq x$, and so   
$E_{\tilde{P}^*}[A_T( \nu^*)]
\leq -x+E_{\tilde{P}^*}[\xi^* +\int_0^T c^*_s ds]
<+\infty,$
which implies that $ A_T( \nu^*) <+\infty\,\, \tilde{P}^* a.s.$ Since $\tilde{P}^*$ is equivalent to $P$, we have $ A_T( \nu^*) <+\infty\,\,P a.s$. This shows that $\nu^*\in \Nc(\tilde{K})$ and so $ \tilde{P}^*\in \bar{\Pc}^{0}$.
By Fatou's Lemma, we have $E[(\tilde{Z}^{*}_T)^\eta]\leq \Liminf_{n\longrightarrow \infty}E[(\bar{Z}^{n}_T)^\eta]<\infty$. Similarly 
$E[(\tilde{Z}^{*}_T)^{1-\bar \eta}]\leq \Liminf_{n\longrightarrow \infty}E[(\bar{Z}^{n}_T)^{1-\bar \eta}]<\infty$. Since $\bar \nu^n\in G_{equi}$ P a.s. and $\bar \nu^n$ converges to $\nu^*$ $dt\otimes dP$ a.e. then  $\nu^*\in G_{equi}$ P a.s. This shows that $\tilde {P}^{*}\in \Pc^0$ and the optimality is deduced.
\ep\\
\\\\
Our aim is to derive a necessary and sufficient condition of optimality of $(c^*,\xi^*)$.
We follow the approach of Duffie and Skiadas \cite{duf1} and El Karoui et al. \cite{elk}, by studying an auxiliary optimization problem without constraints.
Let $\lambda$ be a positive constant, we consider the following  consumption-investment problem
\begin{eqnarray}\label{reward2}
\Sup_{(c,\xi)\in  \Ac }J(x,c,\xi,\tilde P^*,\lambda),
\end{eqnarray}
where the functional $J$ is defined on $ \Ac$ by
\begin{eqnarray} \label{defJ}
J(x,c,\xi,\tilde P^*,\lambda)=Y_0^{x,c,\xi}+\lambda \big(x-E_{\tilde P^*}[\xi+\int_0^Tc_tdt-A_T(\nu^*)]\big).
\end{eqnarray}
We recall the following classical result of convex analysis (see e.g. Luenberger \cite{lue},
Theorem 1 page 217 and Theorem 2 page 221) which relates the solutions of the problems \reff{reward}
and \reff{reward2}.
\begin{Proposition}\label{optim}
We fix $ x\geq v(0,0)$. We assume that the discounting factor is bounded {\bf (H1)}, the set of controls $\Ac$ is non-empty {\bf(H3)}, 
the utility functions satisfy the usual conditions {\bf (H4)}, the absolute value of the utility functions is quasi-concave {\bf (H5)},  and the translation stability on the
set of admissible strategies {\bf (H6)} holds.\\
(i) There exists a positive constant $\lambda^*$ such that
\begin{eqnarray}\label{=vj}
V(x)=\Sup_{(c,\xi)\in \Ac }J(x,c,\xi,\tilde P^*,\lambda^*).
\end{eqnarray}
(ii) The maximum is attained in \reff{reward2} by $(c^*,\xi^*)$. \\
\end{Proposition}
{\bf Proof.}
(i) and (ii): The set $ \Ac$ is convex. The slater condition for
the optimization problem \reff{reward2} holds since
 the strategy $(\tilde c,\tilde \xi)$ defined by $\tilde \xi=\frac{x}{2 }$ and
$\tilde c_t=\frac{x}{2 T}$, $0\leq t\leq T$  is admissible (i.e. $(\tilde c,\tilde \xi)\in \Ac$).
Using Lemma \ref{lemuni}, the value function \reff{=vj} is finite. From
Luenberger \cite{lue}, Theorem 1 page 217, there exists a positive constant $\lambda^*$
such that equality \reff{=vj}  and the assertion (ii) holds.\\
\ep \\
The next result is a dynamic maximum principle. It relates the utility derivatives of the consumption
and the terminal wealth to the density of probability measure which realizes the maximum in the budget constraint and the density of the probability
measure representing the worst case. The proof is technical and is postponed in the Appendix.
\begin{Theorem}\label{maximumprinciple}
We fix $ x\geq v(0,0)$. We assume that the discounting factor is bounded {\bf (H1)}, the set of controls $\Ac$ is not reduced to the null strategy {\bf(H3)}, 
the utility functions satisfy the usual conditions {\bf (H4)}, the absolute value of the utility functions is quasi-concave {\bf (H5)},  and the translation stability on the
set of admissible strategies {\bf (H6)} holds. \\
Let $(c^*,\xi^*)\in \Ac$ be the optimal consumption
and the optimal terminal wealth for \reff{reward2} with $\lambda=\lambda^*$ 
given in Proposition \ref{optim} . Let $(Y^{x,c^*,\xi^*},Z^{x,c^*,\xi^*})$ be the solution for
the BSDE \reff{bsdech}-\reff{bsdechct}. Then the following maximum principle holds:
\begin{eqnarray*}
\bar{\alpha} Z^{*}_T S^{\delta}_T\bar U^{'}(\xi^*)&=&\lambda^* \tilde Z_T^*\,\,dP \,\,a.s.\\
\alpha Z^{*}_t S^{\delta}_t U^{'}(c^*_t)&=&\lambda^* \tilde Z_t^*,\,\,0\leq t\leq T \,\,dt\otimes dP \,\,a.e.\\
\end{eqnarray*}
where $Z_t^{*}=\Ec_t(-\frac{1}{\beta}M^{Y^{*}}),$ such that  $M_t^{Y^{*}}=\int_0^t Z_s^{x,c^*,\xi^*}dW_s,$   $0\leq t\leq T$, $dt\otimes dP$ a.e.
%where $Z^{*}_t=\Ec(-\frac{1}{\beta}\int_0^t Z_s^{x,c^*,\xi^*}dW_s)$,
\end{Theorem}
\begin{Remark}
Theorem \ref{maximumprinciple} provides a characterization 
of the solution of the primal problem in terms of $\tilde Z^*$ the density of $\tilde{P}^{*}$ and $Z^*$ 
the density of the probability measure associated with the worst scenario. 
It is a generalization of the result of Cvitanic and Karatzas \cite{cvi} (Section 12) when $\alpha=\bar \alpha =1$, $Z^*_t=1$, 
$S_t^{\delta}=1$ $dt\otimes dP$ a.e. for all $t\in [0,T]$. They showed that
\begin{eqnarray*}
\bar U^{'}(\xi^*)&=&\lambda^* \tilde Z_T^*\,\,,dP \,\,a.s.\\
U^{'}(c^*_t)&=&\lambda^* \tilde Z_t^*,\,\,0\leq t\leq T \,\,dt\otimes dP \,\,a.e.\\
\end{eqnarray*}
\end{Remark}
\begin{Remark}\label{unicityproba}
From the dynamic programming principle, the unicity of the optimal strategy $(c^*,\xi^*)$ and the unicity of the probability measure $Q^*$ associated with the worst case, we deduce the unicity of the probability measure $\tilde P^*$.
\end{Remark}
\section{Forward-Backward System and Examples}
In this section, we characterize the optimal consumption-investment strategy as the unique solution
of a forward-backward system. This characterization is a consequence of the maximum principle.
In fact, from
Theorem \ref{maximumprinciple}, The optimal terminal wealth $\xi^*$ and the optimal consumption
$c^*_t$ are given by
\begin{eqnarray}
c^*_t&=&I_1\Big(\frac{\lambda^*}{ {\alpha}}S_t^{\delta}
\tilde Z_t^*Z^{*-1}_t\Big),\,\,dt\otimes dP \,\,a.e.\, t\in[0,T]\label{stropc}\\
\xi^*&=&I_2\Big(\frac{\lambda^*}{ \bar{\alpha}} S_T^{\delta}
\tilde Z_T^*Z^{*-1}_T\Big),\,\,dP\,a.s.\label{stropxi}
\end{eqnarray}
where $I_1$  (resp. $I_2$) is the inverse of the derivative function of $U$   (resp.  $\bar U$).
The following result is a direct consequence of Theorem \ref{bd} and Theorem \ref{main1}.
\begin{Theorem}\label{cara}
We fix $ x\geq v(0,0)$. We assume that the discounting factor is bounded {\bf (H1)}, the set of controls $\Ac$ is non-empty {\bf(H3)}, 
the utility functions satisfy the usual conditions {\bf (H4)}, the absolute value of the utility functions is quasi-concave {\bf (H5)},  and the translation stability on the
set of admissible strategies {\bf (H6)} holds. \\
We consider $Y\in \Dc_0^{exp}$,
$Z^{Y}=(Z^{Y}_t)_{t\in [0,T]}$ $\R^d$-valued adapted process satisfying $E[\int_0^T|Z^{Y}_t|^2dt]<\infty$,
$(c^*,\xi^*)\in \Ac(x)$ and $(Z,\tilde{Z})$ two densities of a probability measures equivalent to $P$. 
%and $Z^\nu$ given by (\ref{znu}).
Then, $Y$ coincides with the optimal value process given by $Y^{x,c^*,\xi^*}$, $(c^*,\xi^*)$ are given by \reff{stropc}-\reff{stropxi}, $Z$ coincides with the density
of the minimizing measure $Z^*$ given by $(\ref{minimizingmeasure})$ and  $\tilde{Z}$ coincides with $\tilde Z^*$ given by $(\ref{Ztilde})$, if and only if there
exists $H^*\in \tilde{\Hc}$ satisfying $ X_T^{x,c^*,H^*}\geq \xi^* $ where 
$X_T^{x,c^*,H^*} = x+ \int_0^T H^*_t \mbox{diag}(S_s)^{-1}dS_s -\int_0^T c^*_sds$,
and the following forward-backward system
\begin{eqnarray}\label{densa}
\left\{
\begin{array}{ll}
%dX_t = H_t^{\prime}\mbox{diag}(S_t)^{-1}dS_t -c_tdt,\,\,\,\,\, &X_0^{x,c,H}=x\nonumber\\
dY_t=(\delta_tY_t-\alpha U(c^*_t) +\frac{1}{2\beta}|Z^Y_t|^2)dt+Z^{Y\prime}_tdW_t,&Y_T =\bar \alpha \bar U(\xi^*)\nonumber\\
dZ_t=-\frac{1}{\beta}Z_tZ^{Y\prime}_tdW_t,&Z_0=1\\
d\tilde{Z}_t=\tilde{Z}_t(b_t+\sigma_t^{-1}\nu^*_t)^{\prime}dW_t& \tilde{Z}_0=1
\end{array}
\right.
\end{eqnarray}
admits a unique solution.
\end{Theorem}
{\bf Proof.}
If $Y$ is given by \reff{robust} and $Z$ is given by \reff{minimizingmeasure},
then, from Bordigoni et al. \cite{brom1}, $(Y,Z^Y)$ is the unique solution of
BSDE \reff{BSDE} with terminal condition \reff{BSDEct} and  $Z$ is the solution
of the second equation in our forward-backward system. 
From Theorem \ref{main1}, 
the couple $(c^*,\xi^*)$ is the unique solution of \reff{reward} i.e. $V(x)=Y_0^{x,c^*,\xi^*}$ wich implies that $(Y,Z^Y)$ is the solution of the BSDE
\begin{eqnarray*}
\left\{
\begin{array}{lll}
dY_t&=&(\delta_tY_t-\alpha U(c^*_t) +\frac{1}{2\beta}|Z^Y_t|^2)dt+Z^{Y\prime}_tdW_t\\
Y_T &=&\bar \alpha \bar U(\xi^*),
\end{array}
\right.
\end{eqnarray*}
and there exists $H^*\in \tilde{\Hc}$ such that  
$X_T^{x,c^*,H^*} = x+ \int_0^T H^*_t\mbox{diag}(S_s)^{-1}dS_s -\int_0^T c^*_sds$ and $ X_T^{x,c^*,H^*}\geq \xi^* $. \\
Since $Z$ coincides with $\tilde{Z}^*$, then $Z$ is the solution of the following forward SDE 
\begin{eqnarray*}
\left\{
\begin{array}{lll}
dZ_t=-\frac{1}{\beta}Z_tZ^{Y\prime}_tdW_t\\
Z_0=1.
\end{array}
\right.
\end{eqnarray*}
From Theorem \ref{bd}, $\tilde{Z}$ evolves according to the following forward SDE  
\begin{eqnarray*}
\left\{
\begin{array}{lll}
d\tilde{Z}_t&=&\tilde{Z}_t(b_t+\sigma_t^{-1}\nu^*_t)^{\prime}dW_t\\
\tilde{Z}_0&=&1.
\end{array}
\right.
\end{eqnarray*}
The converse sense is straightforward.
\ep
\begin{Example}
\underline{Incomplete market}: 
In this example, we give an explicit formula for the investment strategy in the risky assets. 
We consider a financial market consisting of two risky assets $S=(S_t^1,S_t^2)_{0\leq t\leq T}$,
where the price is governed by 
\begin{eqnarray*}
\left\{
\begin{array}{ll}
dS^1_t&=b^1 S^1_tdt+\sigma^1 S^1_tdW^1_t,\\
dS^2_t&=b^2 S^2_tdt+\sigma^2 S^2_tdW^2_t,
\end{array}
\right.
\end{eqnarray*}
$(W^1,W^2)$ is a $P$-$\F$ standard Brownian motion, $ b^1, b^2, \sigma^1, \sigma^2$ are constants. 
We take $K=\{h=(h^1, h^2)^{\prime}\in \R^2;\,h_2=0\}$.\\
This is the incomplete market case studied by Karatzas et al \cite{kar2} where the investment is restricted only to the first risky asset. 
It follows that the support function of the convex set $-K$ is given by
\begin{eqnarray*}
\left\{
\begin{array}{ll}
\delta^{supp}(h)&=0\mbox{ if },h_1=0\\
\delta^{supp}(h)&=\infty\mbox{ otherwise },
\end{array}
\right.
\end{eqnarray*} 
and 
\begin{eqnarray*}
\tilde K=\{ h\in \R^2;\,h_1=0\}.
\end{eqnarray*} 
We know that the density of the risk neutral measure is given by 
$\tilde Z_T^*=\Ec_T(-\theta^1 W^1-\int(\theta^2+\frac{\nu^2}{\sigma^2})dW^2)$ where $\theta^i=\frac{b^i}{\sigma^i}$, $i=1,2$
and
by the Girsanov theorem,
\begin{eqnarray*}
\left\{
\begin{array}{ll}
\tilde W_t^1&=W_t^1+\theta^1 t\\
\tilde W_t^2&=W_t^2+\theta^2 t+\int_0^t\frac{\nu^2_s}{\sigma^2}ds
\end{array}
\right.
\end{eqnarray*} 
is a $\tilde P^*$-$\F$ Brownian motion. \\
We fix $x\geq 0$. If $\delta \equiv 0$, $\alpha = 0$, $\bar \alpha =1$ and $\bar U(z)=\log(z)$, then from the recursive
relation \label{recurciverelation}, we obtain
\begin{eqnarray}\label{ReRe}
Y_0^{x,\xi}=-\beta \log E_P\Big[\exp\Big(
- \frac{1}{\beta}\bar U(\xi) \Big)\Big],
\end{eqnarray}
which is a typical example in the dynamic entropic risk measure. We refer to Barrieu and El Karoui \cite{bar} for more details about risk measures.
The stochastic control problem \reff{reward} is related to the problem
\begin{eqnarray}\label{scexa1}
V^{rm}(x):=\Sup_{\xi \in \Xc(x)}E_P\Big[-\exp\Big(
- \frac{1}{\beta}\bar U(\xi) \Big)\Big],
\end{eqnarray}
where $\Xc(x)=\{\xi\geq 0\,,\xi=x+\int_0^T \frac{H_t^1}{S_t^1} dS_t^1,\,\, \frac{H^1}{S^1}\in L(S^1)\mbox{ and }\Sup_{P^{\nu}\in \Pc^0}E_{P^{\nu}}[\xi]\leq x\}$.\\
$\star$ \underline{First step: we solve the stochastic control problem \reff{scexa1}}\\
The utility function $U^{rm}(z)=-\exp\Big(- \frac{1}{\beta}\bar U(z)\Big)$ is strictly concave and increasing. 
It satisfies the Inada condition. From Kramkov and Schachermayer \cite{krasch99}, 
the dual problem admits a solution i.e. there exists a process $\tilde Z^*$ and a constant $z^*$ and the optimal terminal wealth is given by
\begin{eqnarray}\label{wealth}
\xi^*=I^{rm}(z^* \tilde Z_T^*)\,\, a.s.
\end{eqnarray}
where
$I^{rm}(z)=((U^{rm})^{'})^{-1}(z)=(\frac{1}{\beta})^{\frac{\beta}{1+\beta}}z^{-\frac{\beta}{1+\beta}}$. We know by classical results in duality theory, 
that the optimal wealth process of the investor is given by 
\begin{eqnarray*}
X_t^{x,*}=
E_{\tilde P^*}[I^{rm}(z^* \tilde Z_T^*) |\Fc_t].
\end{eqnarray*} 
Since the market is incomplete, the variation of the process $(X_t^{x,*})_t$ is independent of the Brownian motion $W^2$ which implies that  
$\theta^2+\frac{\nu^2_t}{\sigma^2}=0$ $dt\otimes dP$, $t\in [0,T]$. It yields that
\begin{eqnarray} \label{dynamiquerichesseoptimale2}
X_t^{x,*}&=&(\frac{1}{\beta})^{\frac{\beta}{1+\beta}}{z^*}^{-\frac{\beta}{1+\beta}}
\exp(-\frac{\beta}{(1+\beta)^2}\frac{(b^1)^2}{2 (\sigma^1)^2}T) Z_t^{\beta},
\end{eqnarray}
where $Z_t^{\beta}=\Ec_t(\frac{\beta}{1+\beta}\frac{b^1}{\sigma^1}\tilde W^1)$.
Since $X_0^{x,*}=x$, we have $(\frac{1}{\beta})^{\frac{\beta}{1+\beta}}{z^*}^{-\frac{\beta}{1+\beta}}=x$.
From equation \reff{dynamiquerichesseoptimale2} and using It\^o's formula, we have
\begin{eqnarray*}\label{dynamiquerichesseoptimale3}
dX_t^{x,*}=x\exp(-\frac{\beta}{(1+\beta)^2}\frac{(b^1)^2}{2 (\sigma^1)^2}T)
 \frac{\beta}{1+\beta}\frac{b^1}{(\sigma^1)^2}Z_t^{\beta} \frac{d S_t^1}{S_t^1}.
\end{eqnarray*}
Since $dX_t^{x,*}={H^1_t}^*\frac{ d S_t^1}{S_t^1}$, we have by identification that
\begin{eqnarray*}
{H^1}^*_t=x\exp(-\frac{\beta}{(1+\beta)^2}\frac{(b^1)^2}{2(\sigma^1)^2}T)
\frac{\beta}{1+\beta}\frac{b^1 Z_t^{\beta}}{(\sigma^1)^2}, \quad dt\otimes dP\,a.e., \; \forall \, t \in [0,T],
\end{eqnarray*}
and so the number of shares invested in the risky asset $S^1$, which is denoted by $(\hat \theta_t^{\beta})_{t \in [0,T]}$,  is given by
\begin{eqnarray*}
\hat \theta^{\beta}_t=x\exp(-\frac{\beta}{(1+\beta)^2}\frac{(b^1)^2}{2(\sigma^1)^2}T)
\frac{\beta}{1+\beta}\frac{b^1 Z_t^{\beta}}{(\sigma^1)^2 S_t^1}, \quad dt\otimes dP\, a.e., \; \forall \, t \in [0,T].
\end{eqnarray*}
If we send $\beta$ to infinity, we obtain
\begin{eqnarray}\label{str}
\hat \theta^{\infty}_t=x \frac{b^1 }{(\sigma^1)^2 \tilde Z^{*}_t S_t^1}, \quad dt \otimes dP\, a.e., \; \forall \, t \in [0,T].
\end{eqnarray}
$\star$ \underline{Second step: We check $\tilde P^{*}\in \Pc^0$, $\xi^*\in \Ac(x)$ and the integrability conditions on the upper bound \reff{uupper1}-\reff{uupper2}}\\
Since the optimal control  is given by $\nu^*=(\nu_1^*,\nu_2^*)\equiv(0,-b^2)$, we have $\delta^{supp}(\nu^*_t)=0$ $dt\otimes dP$ a.e. $t\in [0,T]$.
A straightforward calculus shows that $E[(\tilde Z^{*}_T)^{\eta}]<\infty$, $E[(\tilde Z^{*}_T)^{1-\bar \eta}]<\infty$ 
and  $\delta^{supp}(\nu^*)$ is equi-integrable with respect to the Lebesgue measure on $[0,T]$. This shows that $\tilde P^{*}\in \Pc^0$. 
Next, we have to compute 
\begin{eqnarray*}
& &E_P[\exp{|\gamma \bar U(\xi^*)|}]\\
&=&E_P\Big[\exp{\Big|\gamma \Big(\ln(x)-\frac{\beta}{(1+\beta)^2}\frac{(b^1)^2}{2(\sigma^1)^2}T
+\frac{\beta}{1+\beta}\frac{b^1}{\sigma^1}\tilde W_T^1-\frac{1}{2}\frac{\beta}{(1+\beta)^2}(\frac{b^1}{\sigma^1})^2 T\Big)\Big|}\Big]<\infty.
\end{eqnarray*}  
This shows that $\xi^*\in \Ac$ and so $\Ac$ in non-empty. Since the inequality $x\geq v(0,0)=0$ is checked, we have $\Ac(x)$ is non empty and $\xi^*\in \Ac(x)$. 
In addition, we have $A_T(\nu^{*})=0$, then $E^{\tilde P^{*}}[A_T(\nu^*)]=0<\infty$ and $\tilde P^{*}$ 
is in the class of probability measures satisfying \reff{uupper2}.\\
$\star$ \underline{Third step: Interpretation of the results}\\
The equation \reff{str} is coherent with the intuition since when $\beta$ goes to infinity, 
we force the penalty term which appears in the dynamic value process (see equation \reff{robust}) to be equal to zero and so
our model of utility maximization under uncertainty converges to a classical utility maximization problem when the underlying model is known.
The optimal strategy of investment in the first risky asset given in \reff{str} corresponds to the solution 
of utility maximization problem in incomplete market when the utility function $\bar U(x)=\log(x)$. 
Such result could be interpreted as a stability result.
In the context of robust maximization problem, the coefficient $\frac{\beta}{\beta+1}\frac{b^1}{\sigma^1}$ could 
be interpreted as a modified relative risk. Also, one could see such coefficient as a change of the level of the volatility. 
The volatility increases from the level $\sigma^1$ to $\frac{\beta+1}{\beta}\sigma^1$. If $\beta$ is close to 0, then the modified relative risk is small enough 
and the number of shares invested in the first risky asset decreases which is consistent with the intuition since we maximize the worst case.  
\end{Example}

\begin{Example}
\underline{Rectangular convex constraints}
%In this example, we give an explicit formula for the investment strategy in the risky assets. 
We consider a financial market consisting of one risky asset where the price is governed by 
\begin{eqnarray*}
%dS^i_t&= S^i_t(b^i_t dt+\sum_{j=1}^d\sigma^{ij} dW^j_t,),\,\,S^i_0=1,\,i=1,\dots,d.
dS_t&=b S_tdt+\sigma S_tdW_t,
%\end{array}
%\right.
\end{eqnarray*}
$W$ is a $P$-$\F$ standard Brownian motion, $ b, \sigma$ are constants. 
We consider the case where $K=[\alpha,\beta]$; $ -\infty<\alpha\leq 0 \leq \beta\leq +\infty$. We know that   $\delta^{supp}(h)=\beta h^{-}-\alpha h^{+}$ and
$\tilde K=\R.$\\
The density of the risk neutral measure is given by 
\begin{equation*}
\tilde Z^*_T=\Ec_T(-\int(\theta +\frac{\nu^*}{\sigma})dW).
\end{equation*}
%by the Girsanov theorem,
%$\tilde W_t^1=W_t^1+\theta^1 t$ is a $\tilde P^*$-$\F$ Brownian motion. \\
%where $(\tilde W^1_t)_{t\in [0,T]}$ is a $\tilde P^*$ Brownian motion.\\  
We fix $x\geq 0$. If $\delta \equiv 0$, $\alpha = 0$, 
$\bar \alpha =1$, $\bar U(z)=\log(z)$, and the recursive relation (\ref{ReRe}) holds,
the stochastic control problem \reff{reward} is related to the problem
\begin{eqnarray}\label{scexa}
V^{rm}(x):=\Sup_{\xi \in \Xc(x)}E_P\Big[-\exp\Big(
- \frac{1}{\beta}\bar U(\xi) \Big)\Big],
\end{eqnarray}
where $\Xc(x)=\{\xi\geq 0\,,\xi=x+\int_0^T \frac{H_t}{S_t} dS_t,\,\, \frac{H}{S}\in L(S)\mbox{ and }\Sup_{P^{\nu}\in \Pc^0}E_{P^{\nu}}[\xi-A_T(\nu)]\leq x\}$.\\
%$\star$ \underline{First step: we solve the stochastic control problem \reff{scexa}}\\
The utility function $U^{rm}(z)=-\exp\Big(-\frac{1}{\beta}\bar U(\xi) \Big)$ is strictly concave and increasing. It satisfies the Inada condition.\\
Following classical arguments of convex duality, see for example Kramkov and Schachermayer \cite{krasch99} and Pham \cite{pha00b}, 
there exists $z^*$, s.t. the optimal wealth process is given by  
\begin{eqnarray*}
\left\{
\begin{array}{ll}
\tilde \xi^*=I^{rm}( z^* \tilde Z_T^*)\,\, a.s.\\
\tilde X_t^{*}=E_{\tilde P^*}[\xi^*-\int_t^T \delta^{supp}(\nu^*_s)ds|\Fc_t]\, dt\otimes dP \mbox{ a.e. }
\end{array}
\right.
\end{eqnarray*} 
The dual problem is given by 
\begin{eqnarray*}
\tilde V^{rm}(z):=\Inf_{\nu \in \Nc(\tilde K)}E_P\Big[\tilde U^{rm}(z Z^{\nu}_T)+z \int_0^T \delta^{supp}(\nu_s)ds\Big],
\end{eqnarray*}
where $\tilde U^{rm}$ is the Fenchel-Legendre transform of $U^{rm}$.\\ 
The dynamic version of the dual control problem is given by 
\begin{eqnarray}\label{DyVer}
\tilde V^{rm}(t,z):=\Inf_{\nu \in \Nc(\tilde K)}E_P\Big[\tilde U^{rm}(z \frac{Z_T^{\nu}}{Z^{\nu}_t})+z \int_t^T \delta^{supp}(\nu_s)ds\Big].
\end{eqnarray}
The HJB equation associated to (\ref{DyVer}) is given by 
\begin{eqnarray}\label{HJBSys}
\left\{
\begin{array}{ll}
-\frac{\partial v }{\partial t}(t,z)+\Sup_{a \in \R}\big[-\Lc^{a}v(t,z)-z(\beta a^{-}-\alpha a^{+})\big]=0,\,\, &(t,z)\in [0,T)\times (0,+\infty)\\
v (T,z)=\tilde U^{rm}(z), \,\, &z \in (0,+\infty).
\end{array}
\right.
\end{eqnarray}
where $$\Lc^{a}v= \frac{1}{2}z^2 \big(\frac{b+a}{\sigma}\big)^2 \frac{\partial^2 v}{\partial^2 z}.$$\\
% -\frac{1}{2}z \frac{b+\nu}{\sigma} \frac{\partial v}{\partial z}$.\\
The HJB equation could be degenerate, the existence of a classical solution is not insured. We should apply the viscosity solutions theory to characterize the dual value function as a viscosity solution of the associated HJB equation. 
\end{Example}
%{\footnotesize
\section{Appendix}
\setcounter{equation}{0}
\setcounter{Assumption}{0}
\setcounter{Example}{0}
\setcounter{Theorem}{0}
\setcounter{Proposition}{0}
\setcounter{Corollary}{0}
\setcounter{Lemma}{0}
\setcounter{Definition}{0}
\setcounter{Remark}{0}
\subsection{ Proof of the comparison theorem}
Denote by
$\Delta Y_t=Y^1_t- Y_t^2$, $\Delta  U_t=\check U^1_t- \check U_t^2$ and $\Delta \bar U_T=\bar U_T^1- \bar U_T^2$.
The pair $(\Delta Y, \int (Z^1-Z^2) dW)$ is the solution of the following equation
\begin{eqnarray*}
d\Delta  Y_t &=&(\delta_t \Delta Y_t -\alpha \Delta U_t)dt
+\frac{1}{2\beta} |Z^1_t|^2 dt
-\frac{1}{2\beta} |Z^2_t|^2dt+(Z^1_t-Z^2_t)dW_t\\
\Delta Y_T&=&\bar \alpha \Delta \bar U_T,\nonumber
\end{eqnarray*}
which implies for any stopping time $T\geq \tau \geq t$, we have
\begin{eqnarray*}
\Delta Y_{\tau}-\Delta Y_{t} &=&\int_t^{\tau}(\delta_s \Delta Y_s -\alpha \Delta U_s)ds
+\frac{1}{2\beta}\int_t^{\tau} |Z_s^1|^2ds\\
&-&\frac{1}{2\beta} \int_t^{\tau}|Z_s^2|^2ds+\int_t^{\tau} (Z^1_s-Z^2_s)dW_s.
\end{eqnarray*}
From the inequality $\int_0^t |Z_s^1|^2ds-\int_0^t|Z_s^2|^2ds-2\int_0^t<Z_s^1,Z_s^2>ds+2\int_0^t|Z_s^2|^2ds=\int_0^t|Z_s^1-Z_s^2|^2ds\geq 0$, 
where $<.,.>$ denotes the inner product associated with the euclidean norm, we deduce that
\begin{eqnarray*}\label{deltav}
\Delta Y_{t} &\leq& \Delta Y_{\tau} -\int_t^{\tau}(\delta_s \Delta Y_s -\alpha \Delta U_s)ds
-\frac{1}{\beta}\int_t^{\tau} <Z^1_s-Z^2_s,Z^{2}_s >ds
-\int_t^{\tau} (Z^1_s-Z^2_s)dW_s.
\end{eqnarray*}
We define the probability measure $ Q^{*,2}$  equivalent to $ P$ where its density is the $P$-martingale
$Z^{*,2}_ T$ with
\begin{eqnarray*}
Z_T^{*,2}=\Ec_T(-\frac{1}{\beta}\int Z_s^2dW_s).
\end{eqnarray*}
Since $\int_0^.(Z_s^1-Z_s^2)dW_s$ is a $P$-martingale, then $ \int_0^.(Z_s^1-Z_s^2)dW_s+\frac{1}{\beta}  \int_0^.<Z_s^1-Z_s^2,Z_s^2>ds$
is a $ Q^{*,2}$-local martingale. Let $(T_n)_n$ be a reducing sequence for \\
$\int_0^.(Z_s^1-Z_s^2)dW_s+\frac{1}{\beta}\int_0^.<Z^1_s-Z^2_s,Z^2>ds$,
then, for $n$ large enough, we have $T_n\geq t$ and so
\begin{eqnarray*}
& &\int_0^t(Z_s^1-Z_s^2)dW_s +\frac{1}{\beta} \int_0^t<Z^1_s-Z^2_s,Z^2_s>ds  \\
&=& E_{Q^{*,2}}[ \int_0^{\tau\wedge T_n}(Z_s^1-Z_s^2)dW_s +\frac{1}{\beta} \int_0^{\tau\wedge T_n}<Z^1_s-Z^2_s,Z^2_s>ds 
|\Fc_t]\,\, \mbox{ on } \{T\geq \tau\wedge T_n \geq t\},
\end{eqnarray*}
which implies
\begin{eqnarray*}
\Delta Y_{t} &\leq& E_{Q^{*,2}}[
\Delta Y_{\tau\wedge T_n } -\int_t^{\tau\wedge T_n }(\delta_s \Delta Y_s -\alpha \Delta U_s)ds
|\Fc_t] \,\, \mbox{ on } \{T\geq \tau\wedge T_n  \geq t\}.
\end{eqnarray*}
Sending $n$ to infinity, we have $\tau\wedge T_n\longrightarrow\tau $, $Q^{*,2}$ a.s. and
$ \Delta Y_{\tau\wedge T_n }\longrightarrow \Delta Y_{\tau}$ $Q^{*,2}$ a.s. Since $Y^1$ (resp. $Y^2$) is in $D_0^{exp}$ and
$ \check U^1 $ (resp. $\check U^2$) is in $D_1^{exp}$, by the dominated convergence theorem, we have
\begin{eqnarray*}
\Delta Y_{t} &\leq& E_{Q^{*,2}}[
\Delta Y_{\tau} -\int_t^{\tau }(\delta_s \Delta Y_s -\alpha \Delta U_s)ds
|\Fc_t] \,\, \mbox{ on } \{T\geq \tau \geq t\}.
\end{eqnarray*}
From the stochastic Gronwall-Bellman inequality ( see Appendix C, Skiadas and Schroder \cite{sch99}), we have
\begin{eqnarray}\label{gbi}
\Delta Y_{t} &\leq&
E_{Q^{*,2}}[
\int_t^{T}\alpha e^{-\int_t^s \delta_s ds}\Delta U_s ds +
\bar \alpha e^{-\int_t^T \delta_s ds} \Delta \bar U_T
|\Fc_t].
\end{eqnarray}
From inequalities \reff{c1}-\reff{c2}, we have
$\Delta Y_{t}\leq 0$, $0\leq t\leq T$,$ \,\,dt\otimes dP \mbox{ a.s. }$ and the result follows.\\
If $Y_0^1=Y_0^2$, then from \reff{gbi}, we have $E_{Q^{*,2}}[
\int_0^{T}\alpha e^{-\int_0^s \delta_s ds}\Delta U_s ds +
\bar \alpha e^{-\int_0^T \delta_s ds} \Delta \bar U_T]=0$, which implies $\Delta U_t=0$, $t\in [0,T]$ and  $\Delta \bar U_T=0$. 
Since the BSDE \reff{BSDE}-\reff{BSDEct} have a unique solution, then $Y_t^1=Y_t^2$, $t\in [0,T]$ a.e.
For the last point, we argue by contradiction. If $Y_0^1=Y_0^2$, then $\Delta U_t=0$, $t\in [0,T]$ and  $\Delta \bar U_T=0$ which contradicts our assumption and so $Y_0^1<Y_0^2$. 
$\hfill \Box$
\subsection{ Proof of the continuity theorem}
We only prove the first statement.
Since $U(c^n_t)\geq U(c_t)$, for all $0\leq t\leq T $ and $\bar U(\xi^n)\geq \bar U(\xi)$,
then from the comparison theorem \ref{ctheorem},
the sequence $((Y^{x,c^n,\xi^n})_{0\leq t\leq T})_n$ is also non-increasing and so
\begin{eqnarray}\label{monoY}
Y^{x,c^1,\xi^1}_t\geq Y^{x,c^n,\xi^n}_t\geq Y^{x, c,\xi}_t,\,\,\, 0\leq t\leq T.
\end{eqnarray}
We define $(Y^{(\infty)}_t)_{0\leq t\leq T}$ as follows:
$Y^{(\infty) }_t=\Lim_{n \longrightarrow \infty}Y^{x,c^n,\xi^n}_t$, $0\leq t\leq T$.
From the definition of $(Y^{x,c^n,\xi^n}_t)_{0\leq t\leq T}$, we have
\begin{eqnarray*}
Y^{x,c^n,\xi^n}_t=-\beta \log E_P[\exp(-\frac{1}{\beta}\int_t^T(\alpha U (c^n_s)
-\delta_sY^{x,c^n,\xi^n}_s)ds)-\frac{1}{\beta}\bar \alpha \bar{U}(\xi^n)|\Fc_t].
\end{eqnarray*}
From inequality \reff{monoY} and the monotonicity property of the sequences $(\xi^n)_n$ and  $\Big((c^n_{t})_{0\leq t\leq T}\Big)_n$, we have
\begin{eqnarray*}
& &\Big|\exp(-\frac{1}{\beta}\int_t^T(\alpha  U(c^n_s)-\delta_s Y^{x,c^n,\xi^n}_s)ds)
-\frac{1}{\beta}\bar \alpha\bar{U} (\xi^n)\Big|\\
&\leq& \exp(\frac{1}{\beta}\int_t^T(\alpha | U(c^n_s)|+\delta_s|Y^{x,c^n,\xi^n}_s|)ds)
+\frac{\bar \alpha}{\beta}|\bar{U} (\xi^n)|)\\
&\leq& \exp\Big(\frac{\alpha}{\beta}\int_0^T(| U(c^1_s)|+|U(c_s)|)ds
+\frac{\|\delta\|_{\infty}T}{\beta}\big(\textrm{ess}\sup_{0\leq t\leq T}|Y^{x,c^1,\xi^1}_t|
+\textrm{ess}\sup_{0\leq t\leq T}|Y_t^{x,c, \xi}|\big)\\
&+&\frac{\bar \alpha}{\beta}(|\bar{U} (\xi^1)|+|\bar{U}(\xi)|)\Big):=g_T.
\end{eqnarray*}
From Cauchy Schwarz inequality, we have
\begin{eqnarray}\label{gT}
E_P[|g_T|] &\leq& E_P\Big[\exp\Big(\frac{2\alpha}{\beta}\int_0^T(| U(c^1_s)|+|U(c_s)|)ds\Big)\Big]^{\frac{1}{2}}\nonumber\\
&& E_P\Big[\exp\Big( \frac{2\|\delta\|_{\infty}T}{\beta}\big(\textrm{ess}\sup_{0\leq t\leq T}|Y^{x,c^1,\xi^1}_t|
+\textrm{ess}\sup_{0\leq t\leq T}|Y_t^{x,c, \xi}|\big)+\frac{2\bar \alpha}{\beta}(|\bar{U} (\xi^1)|+|\bar{U}(\xi)|)\Big)\Big]^{\frac{1}{2}}\nonumber\\
&\leq & E_P\Big[\exp\Big(\frac{2\alpha}{\beta}\int_0^T(| U(c^1_s)|+|U(c_s)|)ds\Big)\Big]^{\frac{1}{2}}\nonumber\\ 
&&E_P\Big[\exp\Big( \frac{4\|\delta\|_{\infty}T}{\beta}\big(\textrm{ess}\sup_{0\leq t\leq T}|Y^{x,c^1,\xi^1}_t|
+\textrm{ess}\sup_{0\leq t\leq T}|Y_t^{x,c, \xi}|\big)\Big)\Big]^{\frac{1}{4}}\nonumber\\
&&E_P\Big[\exp\Big(\frac{4\bar \alpha}{\beta}(|\bar{U} (\xi^1)|+|\bar{U}(\xi)|)\Big)\Big]^{\frac{1}{4}}.
\end{eqnarray}
From the boundedness on the discounting factor {\bf (H1)} and since $(c,\xi)\in \Ac(x)$, $(c^1,\xi^1)\in \Ac(x)$, 
$Y^{x,c, \xi}\in D_0^{exp}$ and $Y^{x,c^1, \xi^1}\in D_0^{exp}$, we have $g_T\in L^1(P)$.
By the dominated convergence theorem, we have
\begin{eqnarray*}
Y^{(\infty)}_t=-\beta \log E_P[\exp(-\frac{1}{\beta}\int_t^T(\alpha U(c_s)-\delta_sY_s^{(\infty)})ds)
-\frac{1}{\beta}\bar{\alpha}\bar{U}(\xi)|\Fc_t],\,\,0\leq t\leq T.
\end{eqnarray*}
Since there exists a unique solution to the BSDE \reff{bsdech}-\reff{bsdechct},
we have necessarily $Y^{(\infty)}=Y^{x,c,\xi}$ and the result follows.
\subsection{Proof of the maximum principle}
We fix $\epsilon>0$ and $\eta>0$ such that $\epsilon <\eta$.\\
{\bf First step:} We prove that
\begin{eqnarray}\label{pm1}
\bar{\alpha} Z^{*}_TS^T_{\delta}\bar U^{'}(\xi^*)&\leq&\lambda^* \tilde Z_T \,\,dP \,\,a.s.
\end{eqnarray}
We consider the following set
\begin{eqnarray*}
A_{\epsilon,\eta} :=\Big \{Z_T^{*} S_{\delta}^T\bar{\alpha}
\bar{U}'(\xi^*)-\lambda^* \tilde Z_T  >0, \epsilon < \xi^* < \eta \Big \}.
\end{eqnarray*}
We define $\xi_n$ as follows: $\xi_n=\xi^* +\frac{1}{n} {\bf 1}_{A_{\epsilon,\eta}}  $.\\
$\star$ \underline {We prove that $(c^*,\xi_n)\in \Ac$:}
From the representation theorem under $\tilde P^*$, there exists a process $H_n\in \tilde \Hc$ such that
\begin{eqnarray*}
\frac{1}{n} {\bf 1}_{A_{\epsilon,\eta}}&=&E_{\tilde P^*}[\frac{1}{n} {\bf 1}_{A_{\epsilon,\eta}}]+\int_0^T{H_n}_s (diag S_s)^{-1}dS_s,
\end{eqnarray*}
which implies that
\begin{eqnarray*}
\xi_n=x+ E_{\tilde P^*}[\frac{1}{n} {\bf 1}_{A_{\epsilon,\eta}}]+ \int_0^T{H_n^*}_s (diag S_s)^{-1}dS_s-\int_0^Tc_s^*ds,
\end{eqnarray*}
where $H_n^*=H^*+H_n$.
%We have
%\begin{eqnarray}\label{maj}
%v(c^*,\xi_n) &=& \sup_{P^{\nu}\in\Pc^0} E^{P^{\nu}} \left[\xi_n+\int_0^Tc^*_tdt -\tilde A_T(P^{\nu}) \right]\nonumber\\
%&=&x+\sup_{P^{\nu}\in\Pc^0} E^{P^{\nu}} \left[ \int_0^T  H_{n,t}dS_t -\tilde A_T(P^{\nu})\right]\nonumber\\
%&\leq&x+ \sup_{P^{\nu}\in\Pc^0} E^{P^{\nu}} \left[ \int_0^T  H_t^*dS_t\right]
%+ \frac{1}{n}\sup_{P^{\nu}\in\Pc^0} E^{P^{\nu}} \left[\int_0^T  \tilde H_t dS_t\right].
%\end{eqnarray}
%Since $(c^*,\xi^*)\in \Ac(x)$ and using Proposition \reff{dualdomin}, we have 
%\begin{eqnarray}\label{stradm}
%\sup_{P^{\nu}\in\Pc^0} E^{P^{\nu}} \left[\xi^*+\int_0^Tc^*_tdt -\tilde A_T(P^{\nu}) \right]\leq x.
%\end{eqnarray}
%From inequality \reff{stradm}, and since we have
%\begin{eqnarray*}
%\sup_{P^{\nu}\in\Pc^0} E^{P^{\nu}} \left[{\bf 1}_{A_{\epsilon,\eta}}
%-\sup_{P^{\nu}\in\Pc^0}E_{\tilde P^{\nu}}[{\bf 1}_{A_{\epsilon,\eta}} ]\right]\leq 0,
%\end{eqnarray*}
%then inequality \reff{maj} implies 
%\begin{eqnarray}\label{stradmi2}
%v(c^*,\xi_n)\leq x.
%\end{eqnarray}
For n large enough, we have $0\leq \frac{1}{n}\leq \frac{\epsilon}{2}$ and so 
\begin{eqnarray*}
{\epsilon}\leq \xi_n\leq \eta+\frac{\epsilon}{2}\mbox{ on the set }\{ \epsilon <\xi^*<\eta\}.
\end{eqnarray*}
From the standard assumptions on the utility functions {\bf (H4)}, we have
\begin{eqnarray*}
\bar U(\epsilon)\leq  \bar U(\xi_n)\leq \bar U(\eta+\frac{\epsilon}{2}) \mbox{ on the set }\{ \epsilon <\xi^*<\eta\}.
\end{eqnarray*}
%and
%\begin{eqnarray*}
%\bar U^{'}(\eta+\frac{\epsilon}{2})\leq \bar U^{'}(\xi_n)\leq \bar U^{'}({\epsilon})\mbox{ on the set }\{ \epsilon <\xi^*<\eta\}.
%\end{eqnarray*}
and so for $n$ large enough, $E[\exp{(\gamma |\bar U(\xi_n)|)}]$ 
%and $E[\exp{(\gamma |\bar U^{'}(\xi_n)|)}]$ are 
is finite, which implies that $(c^*,\xi_n)\in \Ac$.\\
$\star$ \underline{ We prove that $P(A_{\epsilon,\eta})=0$:}
From the definition of $J$ (see \reff{defJ}) and the optimality of the strategy $(c^*,\xi^*)$, we have
\begin{eqnarray} \label{JHn}
0&\geq& n(J(x,c^*,\xi^{n} ,\tilde P^*, \lambda^*)-J(x,c^*,\xi^{*},\tilde P^*,\lambda^*))\\
&=& n(Y_0^{x,c^*,\xi^{n}}- Y_0^{x,c^*, \xi^{*}}) -\lambda^*E_{\tilde P^*}[{\bf 1}_{A_{\epsilon,\eta}}]\nonumber \\
&\geq & nE_{Q^{n}}\Big[ \bar \alpha S_T^\delta
(\bar U(\xi^{n})-\bar U(\xi^*)) \Big]-\lambda^*E_{\tilde P^*}[{\bf 1}_{A_{\epsilon,\eta}}]\nonumber \\
&=& nE_{P}\Big[ Z^{Q^{n}}_T\bar \alpha S_T^\delta
(\bar U(\xi^{n})-\bar U(\xi^*)) {\bf 1}_{A_{\epsilon,\eta}}\Big]-\lambda^*E_{\tilde P^*}[{\bf 1}_{A_{\epsilon,\eta}}],\nonumber
\end{eqnarray}
where the probability measure $Q^{n}$ has a density given by the $P$-martingale
$Z^{Q^{n}}=(Z_t^{ Q^{n}})_{0\leq t\leq T}=(\Ec_t(-\frac{1}{\beta}
M^{x,c^*,\xi^{n}}))_{0\leq t\leq T}$ and $M_t^{x,c^*,\xi^{n}}=\int_0^tZ_s^{x,c^*,\xi^{n}}dW_s$.  \\
Since there exists $\theta^n$ between $\xi^{n}$ and $\xi^*$ such that
$ \bar U(\xi^{n})-\bar U(\xi^*)=\bar U^{'}(\theta^n)(\xi^{n}-\xi^*)$, we deduce that
\begin{eqnarray}\label{cvxi}
n (\bar U(\xi^{n})-\bar U(\xi^*)){\bf 1}_{A_{\epsilon,\eta}}\longrightarrow \bar U^{'}(\xi^*){\bf 1}_{A_{\epsilon,\eta}} \,\,dP \,\,a.s.
\end{eqnarray}
and
\begin{eqnarray}\label{majxi}
| n(\bar U(\xi^{n})-\bar U(\xi^*)){\bf 1}_{\{\epsilon<\xi^*<\eta\}} |\leq \bar U^{'}(\epsilon)\,\,dP \,\,a.s.
\end{eqnarray}
From the definition of $Z_t^{Q^{n}}$, we have
\begin{eqnarray}\label{zqn}
Z_t^{Q^{n}}=\exp(-\frac{1}{\beta} M_t^{x,c^*,\xi^{n}}
-\frac{1}{2\beta^2}<M^{x,c^*,\xi^{n}}>_t).
\end{eqnarray}
From the BSDE \reff{bsdech}, we obtain
\begin{eqnarray}\label{yn}
Y^{x,c^*,\xi^{n}}_t-Y^{x,c^*,\xi^{n}}_0&=&
\int_0^t(\delta_s Y^{x,c^*,\xi^{n}}_s-\alpha U(c_s^{*}))ds
+\frac{1}{2\beta}\langle M^{x,c^*,\xi^{n}}\rangle_t +M_t^{x,c^*,\xi^{n}}.
\end{eqnarray}
Plugging \reff{yn} into \reff{zqn}, we obtain
\begin{eqnarray*}
Z_t^{Q^{n}}= \exp \Big(
\int_0^t\frac{1}{\beta}(\delta_s Y^{x,c^*,\xi^{n}}_s-\alpha  U(c_s^{*}))ds
-\frac{1}{\beta} ( Y^{x,c^*,\xi^{n}}_t -Y^{x,c^*,\xi^{n}}_0)
\Big).
\end{eqnarray*}
From Proposition \ref{monotonieYn} (i), we have
\begin{eqnarray}\label{limZQ}
\Lim_{n \longrightarrow \infty} Z_t^{Q^{n}}
&=& \exp \Big(
\int_0^t\frac{1}{\beta}(\delta_s Y^{x,c^*,\xi^*}_s-\alpha U(c^{*}_s))ds
-\frac{1}{\beta} ( Y^{x,c^*,\xi^*}_t -Y^{x,c^*,\xi^*}_0)
\Big)\nonumber\\
&=& Z_t^*\,\,dt\otimes dP \,\,a.s.
\end{eqnarray}
Under the boundedness on the discounting factor {\bf (H1)} and since $(Y^{x,c^*,\xi^{n}}_t)_{0\leq t\leq T} \in D_0^{exp}$, we have
\begin{eqnarray}\label{majq}
|Z_t^{Q^{n}}|\leq \exp \Big(\frac{T}{\beta} ||\delta||_{\infty}\mbox{ess}\Sup_{0\leq t\leq T}|Y_t^{x,c^*,\xi^{n}}|
+\frac{\alpha}{\beta}\int_0^T |U(c^*_s)|ds+\frac{2}{\beta}\mbox{ess}\Sup_{0\leq t\leq T}|Y_t^{x,c^*,\xi^{n}}|\Big).
\end{eqnarray}
From Proposition \ref{monotonieYn} (i), we have $Y_t^{x,c^*,\xi^{1}}\geq Y_t^{x,c^*,\xi^{n}}\geq Y_t^{x,c^*, \xi^*}$ and so
 \begin{eqnarray}\label{majqp}
\mbox{ess}\Sup_{0\leq t\leq T}|Y_t^{x,c^*,\xi^{n}}|
\leq \mbox{ess}\Sup_{0\leq t\leq T}|Y_t^{x,c^*,\xi^{1}}|
+\mbox{ess}\Sup_{0\leq t\leq T}|Y_t^{x,c^*, \xi^*}|.
\end{eqnarray}
Using the inequalities \reff{majxi}, \reff{majq} and \reff{majqp}, we have
\begin{eqnarray*}
& &| n(\bar U(\xi^{n})-\bar U(\xi^*)){\bf 1}_{\{\epsilon < \xi^* < \eta \}} |
|Z_t^{Q^{n}}|\\
&\leq& \bar U^{'} (\epsilon)  \exp\Big(\frac{T}{\beta} ||\delta||_{\infty}
(\mbox{ess}\Sup_{0\leq t\leq T}|Y_t^{x,c^*,\xi^{1}}|
+\mbox{ess}\Sup_{0\leq t\leq T}|Y_t^{x,c^*,\xi^*}|)\\
&+&\frac{\alpha}{\beta}\int_0^T |U(c_s^*)|ds+\frac{2}{\beta}(\mbox{ess}\Sup_{0\leq t\leq T}|Y_t^{x,c^*,\xi^{1}}|
+\mbox{ess}\Sup_{0\leq t\leq T}|Y_t^{x,c^*,\xi^*}|)\Big):=g_{T}.
\end{eqnarray*}
From Cauchy Schwartz inequality, we have
\begin{eqnarray}\label{gT1}
E_P[|g_{T}|] 
&\leq& \bar U^{'} (\epsilon) E_P\Big[\exp\Big(\frac{2\alpha}{\beta}\int_0^T|U(c_s^*)|ds\Big)\Big]^{\frac{1}{2}}\\
&& E_P\Big[\exp\Big( \frac{2(2+\|\delta\|_{\infty}T)}{\beta}\big(\textrm{ess}\sup_{0\leq t\leq T}|Y^{x,c^*,\xi^{1}}|
+\textrm{ess}\sup_{0\leq t\leq T}|Y_t^{x,c^*, \xi^*}|\big)
\Big)\Big]^{\frac{1}{2}}\nonumber.
\end{eqnarray}
From the boundedness on the discounting factor {\bf (H1)}, and since $(c^*,\xi^*)\in \Ac$, $(c^*,\xi^{1})\in \Ac$, 
$Y^{x,c^*, \xi^*}\in D_0^{exp}$ and $Y^{x,c^*, \xi^{1}}\in D_0^{exp}$, we have $g_{T}\in L^1(P)$.
By the dominated convergence theorem and substituting inequalities \reff{cvxi} and \reff{limZQ} into \reff{JHn},
we have
\begin{eqnarray*}
0 &\geq& \Lim_{n\longrightarrow \infty} E_{Q^{n}}\Big[\bar{\alpha} S_T^\delta
n(\bar U(\xi^*)-\bar U(\xi_n)){\bf 1}_{\{\epsilon < \xi^* < \eta\}} \Big]-\lambda^*E_{\tilde P^*}[{\bf 1}_{A_{\epsilon,\eta}}]\\
& =&E_{Q^{*}} \Big[ \bar{\alpha} S_T^\delta \bar U^{'}(\xi^*){\bf 1}_{A_{\epsilon,\eta}} \Big]-\lambda^*E_{\tilde P^*}[{\bf 1}_{A_{\epsilon,\eta}}].
\end{eqnarray*}
which implies $P(A_{\epsilon,\eta})=0$
for all $0<\epsilon <\eta <\infty$. Sending $\epsilon \longrightarrow 0$ and $\eta \longrightarrow \infty$, we have 
$A_{\epsilon,\eta} \nearrow \Big \{Z_T^{*} S_{\delta}^T\bar{\alpha}
\bar{U}'(\xi^*)-\lambda^* \tilde Z_T  >0 \Big \} $ and so
inequality \reff{pm1} is proved.\\
{\bf Second step:} We prove that
\begin{eqnarray}\label{pm2}
\bar{\alpha} Z^{*}_T S_T^\delta U^{'}(\xi^*)&\geq&\lambda^* \tilde Z_T^{*} \,\,dP \,\,a.s.
\end{eqnarray}
We consider the following set
\begin{eqnarray*}
B_{\epsilon,\eta}:=\Big \{Z_T^{*}\bar{\alpha} S_T^\delta
\bar{U}'(\xi^*)-\lambda^* \tilde Z_T^{*}  <0, \epsilon < \xi^* < \eta \Big \}.
\end{eqnarray*} We define $\xi_n^{'}$ as follows:
$\xi_n^{'}:=\xi^* -\frac{1}{n}{\bf 1}_{B_{\epsilon,\eta}} $.\\
$\star$\underline {We prove that $(c^*,\xi_n^{'})\in \Ac $:}
As in the first step, for n large enough, we have $0\leq \frac{1}{n}\leq \frac{\epsilon}{2}$ and so 
\begin{eqnarray*}
\frac{\epsilon}{2}\leq \xi_n^{'}\leq \eta \mbox{ on the set }\{ \epsilon <\xi^*<\eta\}.
\end{eqnarray*}
From the standard assumptions on the utility functions {\bf (H4)}, we have
\begin{eqnarray*}
\bar U(\frac{\epsilon}{2})\leq  \bar U(\xi_n^{'})\leq \bar U(\eta) \mbox{ on the set }\{ \epsilon <\xi^*<\eta\}.
\end{eqnarray*}
%and
%\begin{eqnarray*}
%\bar U^{'}(\eta)\leq \bar U^{'}(\xi_n^{'})\leq \bar U^{'}(\frac{\epsilon}{2}) \mbox{ on the set }\{ \epsilon <\xi^*<\eta\}.
%\end{eqnarray*}
This shows that for $n$ large enough, $E[\exp{(\gamma |\bar U(\xi_n^{'})|)}]$ 
%and $E[\exp{(\gamma |\bar U^{'}(\xi_n^{'})|)}]$ are 
is finite and so $(c^*,\xi_n^{'})\in  \Ac$.\\
$\star$\underline {We prove that $P(B_{\epsilon,\eta})=0$:}
From the definition of $J$ (see \reff{defJ}) and the optimality of the strategy $(c^*,\xi^*)$, we have
\begin{eqnarray} \label{JHn2}
0&\geq& n(J(x,c^*,\xi^{'n} ,\tilde P^*, \lambda^*)-J(x,c^*,\xi^{*},\tilde P^*,\lambda^*))\\
&=& n(Y_0^{x,c^*,\xi^{'n}}- Y_0^{x,c^*, \xi^{*}}) +\lambda^*E_{\tilde P^*}[{\bf 1}_{B_{\epsilon,\eta}}]\nonumber \\
&\geq & nE_{Q^{'n}}\Big[ \bar \alpha S_T^\delta
(\bar U(\xi^{'n})-\bar U(\xi^*)) \Big]+\lambda^*E_{\tilde P^*}[{\bf 1}_{B_{\epsilon,\eta}}]\nonumber \\
&=& nE_{P}\Big[ Z^{Q^{'n}}_T\bar \alpha S_T^\delta
(\bar U(\xi^{'n})-\bar U(\xi^*)) {\bf 1}_{B_{\epsilon,\eta}}\Big]+\lambda^*E_{\tilde P^*}[{\bf 1}_{B_{\epsilon,\eta}}],\nonumber
\end{eqnarray}
where the probability measure $Q^{'n}$ has a density given by the $P$-martingale
$Z^{Q^{'n}}=(Z_t^{ Q^{'n}})_{0\leq t\leq T}=(\Ec_t(-\frac{1}{\beta}
M^{x,c^*,\xi^{'n}}))_{0\leq t\leq T}$ and $M_t^{x,c^*,\xi^{'n}}=\int_0^tZ_s^{x,c^*,\xi^{'n}}dW_s$.  \\
Since there exists $\theta^n$ between $\xi^{'n}$ and $\xi^*$ such that
$ \bar U(\xi^{'n})-\bar U(\xi^*)=\bar U^{'}(\theta^n)(\xi^{'n}-\xi^*)$, we deduce that
\begin{eqnarray}\label{cvxi2}
n (\bar U(\xi^{'n})-\bar U(\xi^*)) {\bf 1}_{B_{\epsilon,\eta}}\longrightarrow -\bar U^{'}(\xi^*){\bf 1}_{B_{\epsilon,\eta}} \,\,dP \,\,a.s.
\end{eqnarray}
and
\begin{eqnarray}\label{majxi2}
| n(\bar U(\xi^{'n})-\bar U(\xi^*)){\bf 1}_{\{\epsilon <\xi^*<\eta\}} |\leq \bar U^{'}(\epsilon)\,\,dP \,\,a.s.
\end{eqnarray}
From the definition of $Z_t^{Q^{'n}}$, we have
\begin{eqnarray}\label{zqn2}
Z_t^{Q^{'n}}=\exp(-\frac{1}{\beta} M_t^{x,c^*,\xi^{'n}}
-\frac{1}{2\beta^2}<M^{x,c^*,\xi^{'n}}>_t).
\end{eqnarray}
From the BSDE \reff{bsdech}, we obtain
\begin{eqnarray}\label{yn2}
Y^{x,c^*,\xi^{'n}}_t-Y^{x,c^*,\xi^{'n}}_0&=&
\int_0^t(\delta_s Y^{x,c^*,\xi^{'n}}_s-\alpha U(c_s^{*}))ds
+\frac{1}{2\beta}\langle M^{x,c^*,\xi^{'n}}\rangle_t +M_t^{x,c^*,\xi^{'n}}.
\end{eqnarray}
Plugging \reff{yn2} into \reff{zqn2}, we obtain
\begin{eqnarray*}
Z_t^{Q^{'n}}= \exp \Big(
\int_0^t\frac{1}{\beta}(\delta_s Y^{x,c^*,\xi^{'n}}_s-\alpha  U(c_s^{*}))ds
-\frac{1}{\beta} ( Y^{x,c^*,\xi^{'n}}_t -Y^{x,c^*,\xi^{'n}}_0)
\Big).
\end{eqnarray*}
From Proposition \ref{monotonieYn} (ii), we have
\begin{eqnarray}\label{limZQprime}
\Lim_{n \longrightarrow \infty} Z_t^{Q^{'n}}
&=& \exp \Big(
\int_0^t\frac{1}{\beta}(\delta_s Y^{x,c^*,\xi^*}_s-\alpha U(c^{*}_s))ds
-\frac{1}{\beta} ( Y^{x,c^*,\xi^*}_t -Y^{x,c^*,\xi^*}_0)
\Big)\nonumber\\
&=& Z_t^*\,\,dt\otimes dP \,\,a.s.
\end{eqnarray}
Under the boundedness on the discounting factor {\bf (H1)} and since $(Y^{x,c^*,\xi^{'n}}_t)_{0\leq t\leq T} \in D_0^{exp}$, we have
\begin{eqnarray}\label{majqprime}
|Z_t^{Q^{'n}}|\leq \exp \Big(\frac{T}{\beta} ||\delta||_{\infty}\mbox{ess}\Sup_{0\leq t\leq T}|Y_t^{x,c^*,\xi^{'n}}|
+\frac{\alpha}{\beta}\int_0^T |U(c^*_s)|ds+\frac{2}{\beta}\mbox{ess}\Sup_{0\leq t\leq T}|Y_t^{x,c^*,\xi^{'n}}|\Big).
\end{eqnarray}
From Proposition \ref{monotonieYn} (ii), we have $Y_t^{x,c^*,\xi^{'1}}\leq Y_t^{x,c^*,\xi^{'n}}\leq Y_t^{x,c^*, \xi^*}$ and so
 \begin{eqnarray}\label{majqprime2}
\mbox{ess}\Sup_{0\leq t\leq T}|Y_t^{x,c^*,\xi^{'n}}|
\leq \mbox{ess}\Sup_{0\leq t\leq T}|Y_t^{x,c^*,\xi^{'1}}|
+\mbox{ess}\Sup_{0\leq t\leq T}|Y_t^{x,c^*, \xi^*}|.
\end{eqnarray}
Using the inequalities \reff{majxi2}, \reff{majqprime} and \reff{majqprime2}, we have
\begin{eqnarray*}
& &| n(\bar U(\xi^{'n})-\bar U(\xi^*)){\bf 1}_{\{\epsilon < \xi^* < \eta \}} |
|Z_t^{Q^{'n}}|\\
&\leq& \bar U^{'} (\epsilon)  \exp\Big(\frac{T}{\beta} ||\delta||_{\infty}
(\mbox{ess}\Sup_{0\leq t\leq T}|Y_t^{x,c^*,\xi^{'1}}|
+\mbox{ess}\Sup_{0\leq t\leq T}|Y_t^{x,c^*,\xi^*}|)\\
&+&\frac{\alpha}{\beta}\int_0^T |U(c_s^*)|ds+\frac{2}{\beta}(\mbox{ess}\Sup_{0\leq t\leq T}|Y_t^{x,c^*,\xi^{'1}}|
+\mbox{ess}\Sup_{0\leq t\leq T}|Y_t^{x,c^*,\xi^*}|)\Big):=\tilde g_{T}.
\end{eqnarray*}
From Cauchy Schwartz inequality, we have
\begin{eqnarray}\label{gT}
E_P[|\tilde g_{T}|] 
&\leq& \bar U^{'} (\epsilon) E_P\Big[\exp\Big(\frac{2\alpha}{\beta}\int_0^T|U(c_s^*)|ds\Big)\Big]^{\frac{1}{2}}\\
&& E_P\Big[\exp\Big( \frac{2(2+\|\delta\|_{\infty}T)}{\beta}\big(\textrm{ess}\sup_{0\leq t\leq T}|Y^{x,c^*,\xi^{'1}}|
+\textrm{ess}\sup_{0\leq t\leq T}|Y_t^{x,c^*, \xi^*}|\big)
\Big)\Big]^{\frac{1}{2}}\nonumber.
\end{eqnarray}
From the boundedness on the discounting factor {\bf (H1)}, and since $(c^*,\xi^*)\in \Ac$, $(c^*,\xi^{'1})\in \Ac$, 
$Y^{x,c^*, \xi^*}\in D_0^{exp}$ and $Y^{x,c^*, \xi^{'1}}\in D_0^{exp}$, we have $\tilde g_{T}\in L^1(P)$.
By the dominated convergence theorem and substituting inequalities \reff{cvxi2} and \reff{limZQprime} into \reff{JHn2},
we have
\begin{eqnarray*}
0 &\geq& \Lim_{n\longrightarrow \infty} E_{Q^{'n}}\Big[\bar{\alpha} S_T^\delta
n(\bar U(\xi^*)-\bar U(\xi_n^{'})){\bf 1}_{\{\epsilon < \xi^* < \eta\}} \Big] +\lambda^*E_{\tilde P^*}[{\bf 1}_{B_{\epsilon,\eta}}]\\
& =&E_{Q^{*}} \Big[
-\bar{\alpha} S_T^\delta
\bar U^{'}(\xi^*){\bf 1}_{B_{\epsilon,\eta}}
\Big]+\lambda^*E_{\tilde P^*}[{\bf 1}_{B_{\epsilon,\eta}}].
\end{eqnarray*}
which implies $P(B_{\epsilon,\eta})=0$
for all $0<\epsilon <\eta <\infty$. Sending $\epsilon \longrightarrow 0$ and $\eta \longrightarrow \infty$, we obtain
\begin{eqnarray}\label{mp2a}
Z_T^{*}S_T^{\delta}\bar{\alpha}
\bar{U}'(\xi^*)\geq \lambda^* \tilde Z_T^{*}\,\,\mbox{ on the set }\{\xi^* >0\}\,\, dP\,a.s.
\end{eqnarray}
Since the utility function satisfies the Inada conditions (Assumption {\bf (H4)}), we have $P(\xi^* =0)=0$ and so
$P(B_{\epsilon,\eta})=0$ for all $0<\epsilon <\eta <\infty$.\\ 
$\star$ \underline {We prove inequality \reff{pm2}:}
Sending $\epsilon \longrightarrow 0$ and $\eta \longrightarrow \infty$ we have 
$B_{\epsilon,\eta} \nearrow \Big \{Z_T^{*} S_{\delta}^T\bar{\alpha}
\bar{U}'(\xi^*)-\lambda^* \tilde Z_T  <0 \Big \} $ and so
inequality \reff{pm2} is proved.\\
The result follows from \reff{pm1} and \reff{pm2}. The same argument holds for the consumption process.
$\hfill \Box$
\paragraph{Acknowledgements.} We are very grateful to Nicole El Karoui for helpful
comments and fruitful discussions.\\

\end{document}